\documentclass[final]{siamltex}
\usepackage{amsfonts}
\usepackage{amssymb}

\newcommand{\xdownarrow}[1]{%
  {\left\downarrow\vbox to #1{}\right.\kern-\nulldelimiterspace}
}
\usepackage{titlesec}
\titleformat{\section}
  {\normalfont\Large\bfseries\centering}{\thesection}{1em}{}

\title{Rational curves on generic quintic threefolds}

\author{ B. Wang\\
April 4, 2015}

\begin{document}

\maketitle

\begin{abstract}  
Let $X_0$ be a generic quintic threefold in projective space $\mathbf P^4$ over complex numbers and
 $C_0$ be an irreducible rational curve on $X_0$.  Let

$$c_0: \mathbf P^1\to  C_0\subset X_0$$ be its normalization.   In this paper, we show\par
(1) $c_0$ must be an immersion, i.e. the differential $(c_0)_\ast: T_{t}\mathbf P^1 \to T_{c_0(t)} X_0$ is injective 
at each $t\in \mathbf P^1$,\par
(2) the normal bundle of $c_0$ satisfies
$$H^1(N_{c_0/X_0})=0.$$

\end{abstract}

\section{Introduction}

Throughout the paper, we work over complex numbers $\mathbb C$. Unless it is specified, we use Zariski topology.
The word ``generic", which is also called  ``very general" in literature,   is in the sense of countable Zariski topology.

 Let's have a rigorous statement of the result.  
\par

Let $X_0$ be a generic quintic threefold in $\mathbf P^4$ over complex numbers $\mathbb C$. 
Let $c_0: \mathbf P^1\to X_0$ be a  birational map onto its image. So its image
denoted by $C_0$ is a rational curve.
  The regular map $c_0: \mathbf P^1\to \mathbf P^4$ induces a differential map
\begin{equation}\begin{array}{ccc}
(c_0)_\ast|_t: T_{t}\mathbf P^1 &\rightarrow & T_{c_0(t)}X_0\end{array}
\end{equation}
point-wisely. 
This differential map further induces an injective  morphism on the sheaf module, denoted by $(c_0)_\ast$
\begin{equation}\begin{array}{ccc}
(c_0)_\ast: T_{\mathbf P^1 } &\rightarrow & c_0^\ast(T_{X_0}).
\end{array}\end{equation}
(But $(c_0)_\ast$ may not be injective as a vector bundle morphism).

\bigskip

\begin{theorem}  With above set-up, for a generic $X_0$,

(1) $c_0$ must be an immersion , i.e. there exists a normal bundle 
$$N_{c_0/X_0}$$ over $\mathbf P^1$ uniquely determined by $c_0$ such that
$$\begin{array}{ccccccccc}
0 &\rightarrow & T_{\mathbf P^1} &\stackrel{(c_0)_\ast}\rightarrow &
c_0^\ast(T_{X_0}) &\rightarrow & N_{c_0/X_0}  &\rightarrow & 0,
\end{array}$$
is exact, 
\par
(2) the normal bundle satisfies
\begin{equation} H^1(N_{c_0/X_0})=0. \end{equation}
\par

\end{theorem}

\bigskip

\begin{corollary}
Let $X_0\subset \mathbf P^4$ be a generic quintic threefold, and $c_0$ as in theorem 1.1. 
Then  the normal bundle 
$N_{c_0/X_0}$  as in theorem 1.1  has an isomorphism
$$N_{c_0/X_0}\simeq \mathcal O_{\mathbf P^1}(-1)\oplus \mathcal O_{\mathbf P^1}(-1).$$

\end{corollary}
\bigskip

\begin{proof} of corollary 1.2 following from theorem 1.1:  Notice
$$deg(N_{c_0/X_0})=deg(c_0^\ast(T_{X_0}))-deg(T_{\mathbf P^1})=-2.$$

The bundle $N_{c_0/X_0}$ is over $\mathbf P^1$. It is well-known that it can be split into,
\begin{equation} 
N_{c_0/X_0}\simeq \mathcal O_{\mathbf P^1}(k)\oplus \mathcal O_{\mathbf P^1}(-k-2)
\end{equation}
where $k\geq -1$ is an integer.
By Serre duality
\begin{equation}\begin{array}{cc}
H^1(N_{c_0/X_0})&\simeq H^0( ( N_{c_0/X_0})^\ast\otimes \omega_{\mathbf P^1})\\
& \simeq H^0(\mathcal O_{\mathbf P^1}(-2-k)\oplus \mathcal O_{\mathbf P^1}(k)).\end{array}\end{equation}
By theorem 1.1, $H^1(N_{c_0/X_0})=0$. Hence $-1\leq k\leq -1$. 
Therefore $k=-1$. 

\end{proof}

\bigskip

{\bf Remark}. If $c_0$ is not an immersion, $N_{c_0/X_0}$ is only well-defined as a sheaf module of $\mathcal O_{\mathbf P^1}$.  Since the vector bundle's sheaf of regular sections is a sheaf module, theorem 1.1 still holds for the normal sheaf
$N_{c_0/X_0}$. Then the result $H^1(N_{c_0/X_0})=0$ ($H^1$ of the sheaf module) can be extended  to hypersurfaces 
of any types such as Fano, Calabi-Yau, and of general type.  This extension of  $H^1(N_{c_0/X_0})=0$ alone is already enough to help us to understand many problems in these areas. We'll discuss its details  elsewhere.  \par
\bigskip

\subsection{Outline of the proof}\quad\smallskip

The idea of the proof of  theorem 1.1 is to investigate a ``boundary point" of  the space of  rational maps  lying on generic quintic threefolds with a ``compatification" in the essence of GAUGE LINEAR SIGMA MODEL from string theory [4].  This space is called ``a linear model of stable moduli" in [4].
But we do not use string theory. 
Our  technique is the tool in algebraic geometry,  successive blow-ups along subvarieties containing this ``boundary point".   To guide the blow-ups, we choose  an invariant. 
In the following outline of the proof, we'll only state the construction of such an invariant, and skip the details of blow-ups therefore the ``boundary point" and the ``compatification". \bigskip

Let $S=\mathbf P(H^0(\mathcal O_{\mathbf P^4}(5)))$ be the space of all quintics. 
Theorem 1.1 is stated in terms of rational curves and quintic threefolds in projective space. 
But in this paper we'll stick with the affine space for the simplicity. So let
$$\mathbb C^{5d+5}$$ be the
vector space, 
$$( H^0(\mathcal O_{\mathbf P^1}(d))^{\oplus 5}$$
whose open subset parametrizes  the set of non-constant regular maps $$\mathbf P^1\to \mathbf P^4$$ whose push-forward cycles have degree $d$.\footnote {The automorphism of $\mathbf P^1$ induces a
$PGL(2)$ group action on $\mathbf P(\mathbb C^{5d+5})$.  Let $$PGL(2)( c_0)\subset \mathbf P(\mathbb C^{5d+5})$$ be the orbit of $c_0\in \mathbf P(\mathbb C^{5d+5})$. } 
Throughout the paper, we let $$M= \mathbb C^{5d+5}.$$ 
Let $M_d$ be the subset that consists of all generically one-to-one (to its image) maps $c$ whose images $c_\ast(\mathbf P^1)$ have degree $d$. 
We'll use the affine coordinates of $M$. \footnote {This space $M$ added with infinity  is the compact space we use for our compactification mentioned above.  Our ``boundary point" is not at the infinity and the argument is local. Therefore such a compatification as a global space does not play a role in the proof, and
it will not be mentioned. 
However the projectivization
$\mathbf P(M)$ is the compact space used in GAUGE LINEAR SIGMA MODEL in string theory [4].  The significance of such a specialization is  the same as in GAUGE LINEAR SIGMA MODEL:  $M$ does not parametrize irreducible rational curves of degree $d$. In fact $M$ has a stratification:
$$M\supset \{c: deg(c(\mathbf P^1))\leq d-1 \}\supset \cdots\supset \{c: c(\mathbf P^1)=a\ point\}.$$
Our final ``boundary point" lies in a lower or the lowest stratum. This is one of main difficulties resolved by blow-ups in the proof. }  
 Assume $c_0^\ast(f_0)=0$ for $c_0\in M_{d}$ as in theorem 1.1.
Let $\mathbb L\subset S$ be an open set of the plane spanned by quintics $f_0, f_1, f_2$, where $f_0, f_1, f_2$ are generic in $S$.
Let  $$\Gamma_{\mathbb L}\ni (c_0, [f_0])$$ be an irreducible component of the incidence scheme
\begin{equation}\{(c, [f])\subset M_d\times \mathbb L: c^\ast(f)=0\}
\end{equation} 
that is onto $\mathbb L$, where $[f_0]$ denotes the image of $f_0$ under the map
$$\begin{array}{ccc}
H^0(\mathcal O_{\mathbf P^4}(5))-\{0\} &\rightarrow & S. \end{array}$$
We assume $\Gamma_{\mathbb L}$ exists. Let $P$ be the projection $$\Gamma_{\mathbb L}\to M.$$
The idea of the proof is to show that the scheme,
$$P(\Gamma_{\mathbb L})\subset M$$ is
a reduced, irreducible quasi-affine scheme of dimension $6$.  The method is straightforward to show its defining polynomials at a generic point have
non-degenerate Jacobian matrix (by that we mean it has full rank).\footnote{ In algebra, this actually is a study of the fitting ideal of the module of a
differential sheaf through successive blow-ups. The goal is to show this fitting ideal $Fitt_i$ for the priorly determined set of polynomials has length $i=0$.  }  See
definition 1.10 below for the precise definition of a Jacobian matrix.   All differentials and partial derivatives used throughout the paper are in algebraic sense, i.e.
defined as in [7] (because all functions are holomorphic).  In the following we describe its defining polynomials and a differential form representing
the Jacobian  matrix. 

Choose generic $5d+1$ distinct points  $t_i\in \mathbf P^1(\mathbb C^2)$ (generic in $Sym^{5d+1}(\mathbf P^1(\mathbb C^2))$). Throughout the paper, unless specified otherwise, we'll use $t_i$ to denote a  complex number  which is a point in $\mathbb C\subset \mathbf P^1(\mathbb C^2)$. 
Next we consider differential 1-forms $\phi_i$ on $M$:
  
\begin{equation} \phi_i= d
\left|  \begin{array}{ccc} f_2(c(t_i)) & f_1(c(t_i)) & f_0(c(t_i))\\
f_2(c(t_1)) & f_1(c(t_1)) & f_0(c(t_1))\\
f_2(c(t_2)) & f_1(c(t_2)) & f_0(c(t_2))
\end{array}\right|\end{equation}
for $i=3, \cdots, 5d+1$, and variable $c\in M$, where $|\cdot |$ denotes the determinant of a 
matrix.  Notice $\phi_i$ are uniquely defined provided  the quintics $f_i$ are in an affine open set of $S$ ( $\phi_i$ is not invariant under the $GL(2)$ action of $\mathbb C^2$ ).  Let \begin{equation} \omega(M, \mathbf t)=\wedge_{i=3}^{5d+1}\phi_i \in H^0(\Omega^{5d-1}_M) \end{equation}
be the $5d-1$-form. This $\omega(M, \mathbf t)$  is the dual expression of the Jacobian matrix of some defining polynomials for the scheme
 $P(\Gamma_{\mathbb L})$, where $\mathbf t=(t_1, \cdots, t_{5d+1})$.

So
the crucial point of this definition is that the polynomials inside of ``d" operator,
 \begin{equation} 
\left|  \begin{array}{ccc} 
f_2(c(t_i)) & f_1(c(t_i)) & f_0(c(t_i))\\
f_2(c(t_1)) & f_1(c(t_1)) & f_0(c(t_1))\\
f_2(c(t_2)) & f_1(c(t_2)) & f_0(c(t_2))
\end{array}\right| \end{equation}
 $i=3, \cdots, 5d+1$ generically define the scheme (the scheme-theoretical image under $P$), 
$$P(\Gamma_{\mathbb L}).$$

Using successive blow-ups at a ``boundary point" of $P(\Gamma_{\mathbb L})$, we proved that 
\bigskip

\begin{proposition} For generic choices of quintics, $f_0, f_1, f_2$, and 
$$(t_1, \cdots, t_{5d+1})\in Sym^{5d+1}(\mathbf P^1), $$
$\omega(M, \mathbf t)$ is a non-zero differential form of degree $5d-1$ when restricted to a non empty open set of  $P(\Gamma_{\mathbb L}) $, i.e.
the set of global sections $\{\phi_i\}_{i=3, \cdots, 5d+1}$ is linearly independent in the
$\mathcal O(P(\Gamma_{\mathbb L})  )$ module, 
$$H^0(\Omega_M\otimes \mathcal O_{P(\Gamma_{\mathbb L}) }).$$

\end{proposition}

\bigskip

Next by an argument on Zariski tangent spaces, mainly from the fact that the ideal of scheme,
$$P(\Gamma_{\mathbb L})$$ is generated by 
polynomials,
 \begin{equation} 
\left|  \begin{array}{ccc} f_2(c(t_1)) & f_1(c(t_1)) & f_0(c(t_1))\\
f_2(c(t_2)) & f_1(c(t_2)) & f_0(c(t_2))\\
f_2(c(t_i)) & f_1(c(t_i)) & f_0(c(t_i))\end{array}\right|\end{equation}
we obtain 
\bigskip

\begin{proposition}
If $\omega(M, \mathbf t)$ is non-zero  on $P(\Gamma_{\mathbb L})$(The algebraic equivalence of this is that the set $\{\phi_i\}_{i=3, \cdots, 5d+1}$ 
is linearly independent  in the
$\mathcal O(P(\Gamma_{\mathbb L}))$ module), then the Zariski tangent space of 
$P(\Gamma_{\mathbb L})$ at a generic maximal point must be 
\begin{equation}
dim(M)-deg(\omega(M, \mathbf t))
\end{equation}
\end{proposition}

\bigskip

{\bf Remark} The form $\omega(M, \mathbf t)$ is not invariant under the $GL(2)$  action and depends on $M, \mathbf t, \mathbb L$, but
the zero locus  $\{\omega(M, \mathbf t)=0\}\subset \mathbb C^{5d+5}$ is $GL(2)$ invariant, independent of $M$ (up to an isomorphism) and $\mathbf t$.
\footnote{ The ideal of $\{\omega(M, \mathbf t)=0\}$ is a Jacobian ideal.} The proof of proposition 1.3 is the main body of the paper. It is achieved by successive blow-ups at a rational curve 
lying on  a product of 5 distinct planes in general positions. \bigskip

Continuing from this proposition, by the surjectivity of $\Gamma_{\mathbb L}$ to $\mathbb L$, the argument on Zariski tangent 
spaces shows that the dimension of Zariski tangent space of
$\Gamma_{\mathbb L}$ at a generic point must be the same as that of Zariski tangent space of
$P(\Gamma_{\mathbb L})$ at a generic point. Thus the proposition 1.4 implies that
the dimension of the Zariski tangent space of $\Gamma_{\mathbb L}$ at a generic point is
$6$. 
By lemmas 3.2, 3.3,
this directly leads to theorem 1.1:\par
(1) $c_0$ is an immersion, \par
(2) and
$$H^1(N_{c_0/X_0})=0.$$
\bigskip

\bigskip

The major part of the proof is proposition 1.3. Let's give a very rough sketch of the proof.
In the case $P(\Gamma_{\mathbb L})$ is smooth at $c_2$, the ideal of the proof is simply to take higher partial 
derivatives of $\omega(M, \mathbf t)$ at $c_2$ along directions in $P(\Gamma_{\mathbb L})$ to discover some
higher partial derivative is non zero. This will show $\omega(M, \mathbf t)$ restricted to $P(\Gamma_{\mathbb L})$ is not identically zero. 
In this paper $P(\Gamma_{\mathbb L})$ is most likely non smooth at $c_2$, so we replace partial derivatives by blow-ups.
It is quite easy to see that the centers of the blow-ups must be those ``bad" points making $\omega_(M, \mathbf t)$ vanish, i.e.  the directions of the partial derivatives in $P(\Gamma_{\mathbb L})$ are transversal to the collection of ``bad points".  But the difficulty lies in the fact that
it is not feasible to  take partial derivatives of each entry of the Jacobian.  The invariant to control the succesive blow-ups is  the blow-up element $c_2^{(\kappa)}$ of the rational map $c_2$.  This $c_2^{(\kappa)}$ , 
in case $P(\Gamma_{\mathbb L})$ is smooth at $c_2$,  is exactly the higher partial derivative of family of rational maps $c\in P(\Gamma_{\mathbb L})$ evaluated at $c_2$. Then it is easy to see that this $c_2^{(\kappa)}$ is not determined by $c_2$ which could be completely trivial (constant map). It is determined by the family, $ P(\Gamma_{\mathbb L})$ (the lowest desired order of deformation of $c_2$ in this family). 
The lemma 5.3 proves that the blow-up of $\omega(M, \mathbf t)$ (the Jacobian) evetually becomes nonzero when
the partial derivative $$c_2^{(\kappa)}=((  c_2^0)^{(\kappa)}, \cdots,  ( c_2^4)^{(\kappa)})$$ has
\par
(a)  all  components $(c_2^i)^{(\kappa)}$  are  non-zero, \par
(b)  $5d$ zeros of $(c_2^i)^{(\kappa)}(t)=0$ are distinct.

\par
In step 1 of lemma 5.3, we reach part (a). In step 2 of lemma 5.4, we arrive at (b).

\bigskip

To summarize it, in this paper we investigate two distinct invariants: incidence scheme $\Gamma$ and normal bundle
$N_{c_0/X_0}$.  Section 3 gives a connection between these two. But almost entire paper concentrates  on a result of $\Gamma$ only. 
The theorem 1.1 is just the expression of this result on $\Gamma$ in terms of the normal bundle $N_{c_0/X_0}$.
 \bigskip

Also our result just provides some bases for Gromov-Witten invariants, but it does not address  them.

\bigskip

\subsection{Technical notations}\quad\smallskip

In this section, we collect all technical notations and definitions used in this paper. Some of them may already be defined before.\bigskip

{\bf Notations}:
\par
(1) $S$ denotes the space all quintics, i.e. $S=\mathbf P(H^0(\mathcal O_{\mathbf P^4}(5)))$.\par
Let $[f]$ denote the image of $f$ under the map
$$\begin{array}{ccc}
H^0(\mathcal O_{\mathbf P^4}(5))-\{0\} &\rightarrow & S. \end{array}$$

(2) Let $$M$$ be
$$\mathbb C^{5d+5}\simeq (H^0(\mathcal O_{\mathbf P^1}(d))^{\oplus 5}$$ 
and $M_{d}$ be the subset that parametrizes  all birational-to-its-image maps $$\mathbf P^1\to \mathbf P^4$$ whose push-forward cycles  have degree $d$. 
 \par
(3) Throughout the paper, if $$c:  \mathbf P^1\to \mathbf P^4, $$
is regular,   $c^\ast(\sigma)$   denotes the pull-back section  of
section $\sigma$ of some bundle over $\mathbf P^4$. The vector bundles will not always be specified, but they are apparent in the context.\par
(4) Let $Y$ be a scheme,  $y\in Y$ be a closed point, $Z\subset Y$ be a subscheme (open or closed) and $\mathcal M$ be a quasi-coherent
sheaf of $\mathcal O_{Y}$-module.
Then $\mathcal O_{y, Y}$ denotes the local ring, $\Omega_{Y}$ denotes the sheaf of differentials, 
  $\mathcal M|_{(Z)}$ denotes the inverse sheaf module $i^\ast(\mathcal M)$ where $i: Z\hookrightarrow Y$ is the embedding. We call  
$\mathcal M|_{(Z)}$ the restriction of $\mathcal M$ to $Z$. 
$\mathcal M|_{Z}$ denotes the localization of $\mathcal M$ at $Z$, which is a $\mathcal O_{Z, Y}$ module.  Thus
$$\mathcal M|_{(\{y\})}=\mathcal M|_Z\otimes k(y),$$
where $k(y)$ is the residue field of the maximal point $\{y\}$.  
\par
If $Y$ is quasi-affine scheme, $\mathcal O (Y)$ denotes the ring of regular functions on $Y$.  

\par

(5) Let $\alpha\in T_{c_0}M$, and $$g: M\to H^0(\mathcal O_{\mathbf P^1}(r))$$
be a regular map. 
Then the  image $g_\ast(\alpha)$ of $\alpha$ under the differential map at $c_0$ is denoted by
\begin{equation}
{\partial g(c_0(t))\over \partial \alpha}\in T_{g(c_0)}(H^0(\mathcal O_{\mathbf P^1}(r)))=H^0(\mathcal O_{\mathbf P^1}(r)).
\end{equation}
(use the identification $T_{g(c_0)}(H^0(\mathcal O_{\mathbf P^1}(r)))=T_0(H^0(\mathcal O_{\mathbf P^1}(r)))$).
\par
(6) If $Y$ is a scheme, $|Y|$ denotes the induced reduced scheme of $Y$.

\bigskip

\begin{definition} By the adjunction formula,  $deg(N_{c_0/X_0})=-2$. So 
let
\begin{equation}
N_{c_0/X_0}\simeq \mathcal O_{\mathbf P^1}(k)\oplus  O_{\mathbf P^1}(-k-2).\end{equation}
where $k\geq -1$ is an integer determined by the normal bundle $N_{c_0/X_0}$.
We'll fix an isomorphism in (1.13)  throughout. 
Let $E$ be the inverse of the summand 
$\mathcal O_{\mathbf P^1}(k)$ under the map
$$c_0^\ast(T_{X_0})\to N_{c_0/X_0}.$$
Thus  $E$ is a rank $2$ subbundle of $c_0^\ast(T_{X_0})$. 
\bigskip
\end{definition}

\bigskip

So  if $k\geq 0$, $E$ is the sub-bundle of the bundle $c_0^\ast(T_{X_0})$,  generated by all the holomorphic sections of
 $c_0^\ast(T_{X_0})$.
\bigskip

\begin{definition} 
(a) 
 If $f\in H^0(\mathcal O_{\mathbf P^4}(5))$ is a quintic polynomial other than $f_0$, we denote the direction of
the line through two points $[f], [f_0]$ in the projective space, $\mathbf P(H^0(\mathcal O_{\mathbf P^4}(5))$  by $\overrightarrow f$.
So $$\overrightarrow f\in T_{[f_0]}\mathbf P(H^0(\mathcal O_{\mathbf P^4}(5)).$$\par
(b) 
Note that the vector $\overrightarrow f$ is well-defined up-to a non-zero multiple. 
In case when $c_0$ can deform to all quintics to the first order, i.e. the map in (3.1) is surjective, this naturally gives a 
section $<\overrightarrow f>$ of the  bundle $c_0^\ast(T_{\mathbf P^4})$ (may not be unique),  to each deformation $ \overrightarrow f$ of the
quintic $f_0$. This is easily can be understood as the direction of the moving $c_0$ in the deformation $(\overrightarrow f, <\overrightarrow f>)$ of
the pair $(c_0, f_0)$.

\end{definition}

\bigskip

\begin{definition}
Let $\Gamma$  be an irreducible component of the incidence scheme 
\begin{equation}\{(c, f)\subset M_d\times \mathbf P(H^0(\mathcal O_{\mathbf P^4}(5))): c^\ast(f)=0\}
\end{equation} 
that dominates $S=\mathbf P(H^0(\mathcal O_{\mathbf P^4}(5)))$. \end{definition}
Let $(c_0, [f_0])\in \Gamma$ be a generic point.  Throughout the paper we assume that such a $\Gamma$ exists. 
\bigskip

{\bf Remark}: The existence of such a $\Gamma$ is equivalent to  the assumption of theorem 1.1: 
$X_0$ is generic. Results in section 3 only need a weaker assumption, but the main propositions 1.3, 1.4 rely on this stronger assumption---$\Gamma$ exists.  In this paper, to avoid the distraction, we use the unified and consistent assumption--$\Gamma$ exists. 
\bigskip

\begin{definition} Let $f_1, f_2 \in H^0(\mathcal O_{\mathbf P^4}(5))$ be two quintics different from $f_0$.
Let $\mathbb L$ be an open set of the plane in 
$$ \mathbf P(H^0(\mathcal O_{\mathbf P^4}(5)))$$ spanned by $[f_0], [f_1], [f_2]$ and centered around $[f_0]$. 
 \end{definition}

\bigskip

\begin{definition} Let
\begin{equation} 
\Gamma_{\mathbb L}=\Gamma\cap ( M\times \mathbb L)
\end{equation} be an irreducible component of the restriction of $\Gamma$ to $M\times \mathbb L$ such that it is onto ${\mathbb L}$, and
\begin{equation} 
\Gamma_{f_0}, \ for\ generic\ f_0\in \mathbb L
\end{equation} is an irreducible component of
$$P(\Gamma\cap ( M\times \{[f_0]\}))$$
where $P$ is the projection to $M$.

 \end{definition}

\bigskip

\begin{definition} 
Let $V$ be a smooth analytic variety with analytic coordinates $x_1, \cdots, x_n$, Let $f_1, \cdots, f_m$ be holomorphic functions on $V$.
We define

\begin{equation} \begin{array}{c} 
\left (\begin{array}{cccccc}
{\partial f_1\over \partial x_1}  &{\partial f_1\over \partial x_2}   & \cdots & {\partial f_1\over \partial x_n}   \\
{\partial f_2\over \partial x_1}  &{\partial f_2\over \partial x_2} & \cdots & {\partial f_2\over \partial x_n}  \\
 \vdots & \vdots & \cdots &\vdots \\
{\partial f_m\over \partial x_1}  &{\partial f_m\over \partial x_2} & \cdots & {\partial f_m\over \partial x_n}
\end{array}\right). \end{array}\end{equation}

to be the Jacobian matrix of functions $f_1, \cdots, f_m$. This Jacobian matrix depends on the coordinates 
$x_1, \cdots, x_n$. 

\end{definition}

\bigskip
This definition is crucial. It must be noted that the Jacobian matrices depend on the coordinates $x_1, \cdots, x_n$. One of main difficulties of this paper is to search for such coordinates that would make Jacobian matrices simpler.  One may wish to compare this definition with that of Jacobian ideals which is more intrinsic, but very elusive in our situation.   

\bigskip

In section 2, we prove that original Clemens' conjecture follows from theorem 1.1. In the rest of the paper, we are  going to concentrate on  proving theorem 1.1.  In section 3, we express a first order condition, and investigate the Zariski tangent spaces of incidence schemes. In section 4, we prepare the analytic
coordinates of $M$ for the computation. In section 5, we use
the sheaf of differentials to show the non-vanishing property of $5d$-$1$-form $\omega(M, \mathbf t)$  on the scheme $\mathbf P(\Gamma_{\mathbb L})$. This is the central section of the paper. 
It leads the proof of propositions 1.3, 1.4. Section 6 collects two known examples which emphasize on the  singular rational curves. \bigskip

{\bf Acknowledgment}. We would like to thank Bruno Harris who clears our understanding of the map $(c_0)_\ast$ in (1.2).

\bigskip

\bigskip
\section{Clemens' conjecture}
\quad

Theorem 1.1 and corollary 1.2 prove the Clemens' conjecture that  was first proposed 30 years ago. During the last thirty years, there are many articles on the conjecture. The most of them followed the early idea of Katz ([6]) to show that there is only one irreducible component of the incidence scheme, containing a smooth rational curve and  dominating the space of quintics.       \bigskip

\subsection{Historical remark}\quad\smallskip

Rational curves on hypersurfaces have been great interests for many years in algebraic geometry. The Clemens' conjecture sits in the center
of many major problems in this area. The theorem and the corollary in this paper, stated above are meant to give solutions to the main part of
the conjecture. In order to distinguished the proved and non-proved parts of conjecture, in the following we state the Clemens' original conjecture.  In [2], its original 1986 statement, Clemens 
proposed:\par``(1) the generic quintic threefold $V$ admits only finitely many rational curves of each degree.\par
(2) Each rational curve is a smoothly embedded $\mathbf P^1$ with normal bundle
\begin{equation}
\mathcal O_{\mathbf P^1}(-1)\oplus  O_{\mathbf P^1}(-1).\end{equation}\par
(3) All the rational curves on $V$ are mutually disjoint. The number of rational curves of degree $d$ on $V$ is
\begin{equation}
(interesting \ number)\cdot 5^3\cdot d. "
\end{equation}
\smallskip

In 1995, Vainsencher  found the degree 5, 6-nodal rational curves in the generic quintic threefolds ([8]). This partially disproved part (2) in the Clemens' conjecture and leave
the part (1) unanswered. At the meantime Mirror symmetry came to the stage to redefine the approach in part (3). 
Based on Vainsencher's result, in 1999, motivated by the Gromov-Witten invariants in the mirror symmetry, 
Cox and Katz  modified the Clemens' original conjecture to the most current form ([4]):\par
`` Let $V\subset \mathbf P^4$ be a generic quintic threefold. Then for each degree $d\geq 1$, we have\par

(i) There are only finitely many irreducible rational curves $C\subset V$ of degree $d$.\par
(ii) These curves, as we vary over all degree, are disjoint from each other.\par
(iii) If $c: \mathbf P^1\to C$ is the normalization of an irreducible rational curve $C$, then the normal bundle
has isomorphism
$$N_{c/V}\simeq \mathcal O_{\mathbf P^1}(-1)\oplus \mathcal O_{\mathbf P^1}(-1). "$$

\bigskip

{\bf Remark}. Cox and Katz's conjecture (iii) should be understood as in two steps. First  $N_{c/V}$ must be a locally free sheaf, secondly
$$N_{c/V}\simeq \mathcal O_{\mathbf P^1}(-1)\oplus \mathcal O_{\mathbf P^1}(-1). $$
We proved the first by showing that $c_0$ is an immersion. \par
The conjecture is proved to be correct for $d\leq 9$ by the work of Katz ([6]), Johnsen and Kleiman ([5]), and Cox and Katz ([4]),  etc.

\subsection{A proof of Clemens' conjecture}\quad\smallskip

\bigskip

\begin{corollary}
Let $X_0\subset \mathbf P^4$ be a generic quintic threefold. Then for each degree $d\geq 1$, we have\par
(i) there are only finitely many irreducible rational curves $C_0\subset X_0$ of degree $d$. 
\par
(ii) Each rational curve in (i) is an immersed rational curve with normal bundle
$$N_{c_0/X_0}\simeq \mathcal O_{\mathbf P^1}(-1)\oplus \mathcal O_{\mathbf P^1}(-1).
$$
By  ``immersed rational curve" we mean that the normalization map is an immersion. 
\end{corollary}
\bigskip

\begin{proof} of corollary 2.1 following from theorem 1.1 and corollary 1.2:  The existence of rational curves on a generic quintic was proved in [3], [6].
So it suffices to prove the finiteness. 
Part (i) follows from part (ii). So let's prove part (ii).  Let $C_0$ be an irreducible rational curve of degree $d$ on $X_0$. Then we take a normalization of $C_0$, and denote it by $c_0:\mathbf P^1\to X_0$. Since $X_0$ is generic, we have the set-up for corollary 1.2.  Applying corollary 1.2, we obtain part (ii).    
\end{proof}

 \bigskip

Corollary 2.1  proves the modified Clemens' conjecture, namely parts (i) and (ii)
of Cox and Katz's statements.  Clemens' original conjecture must be modified in the light  of  Vainsencher's result. 

\bigskip

\section{ First order }

\bigskip

\subsection{First order deformations of the pair} \quad\smallskip

 Let's start the problem in its first order.\bigskip

\begin{lemma} Let $f_0$ be a generic quintic threefold containing a rational map $c_0$ as before. 
If $(c_0, [f_0])\in |\Gamma|$ is generic , then the projection
\begin{equation} \begin{array}{ccc}T_{(c_0, [f_0])}\Gamma &\stackrel{P^s}
\rightarrow & T_{[f_0]}S\end{array}\end{equation}
is surjective, where $S=\mathbf P(H^0(\mathcal O_{\mathbf P^4}(5))$.

\end{lemma}

\bigskip

\begin{proof}
Let $|\Gamma |\subset \Gamma$ be the reduced scheme of the scheme $\Gamma$.
By the genericity of $f_0$,  the projection 
\begin{equation}\begin{array} {ccc} |\Gamma | &\rightarrow &
S\end{array}\end{equation}
is dominant. Hence in a neighborhood of a generic point 
$(c_0, [f_0])\in |\Gamma|$, the projection is a smooth map. Thus 
\begin{equation}\begin{array} {ccc} T_{(c_0, [f_0])}|\Gamma | &\rightarrow &
T_{[f_0]}S\end{array}\end{equation}
is surjective. This proves the lemma

\end{proof}

 To elaborate definition 1.6  in introduction, we apply this lemma to obtain that for any $\alpha\in T_{[f_0]}S$,  there is  
a section denoted by 
$$<\alpha>\in H^0(c_0^\ast(T_{\mathbf P^4}))$$ such that $$ (\alpha, <\alpha>)$$ is tangent to the 
universal hypersurface $$\mathcal X=\{(x, f): x\in div(f)\}$$ in $$\mathbf P^4\times S.$$
 Note that $<\alpha>$ is  unique up to a section in $H^0(c_0^\ast(T_{X_0}))$.  But we will always fix $<\alpha>$ as in introduction.
 
\bigskip

\bigskip
\subsection{The incidence scheme}\quad\smallskip

In this subsection, we study the Zariski tangent spaces of various incidence schemes to reveal
a connection between the incidence scheme and the normal  sheaf.  \bigskip

\begin{lemma} Let $[f_0]\in S$ be a generic point, $\mathbb L_1\subset S$ an open set of the pencil containing 
$f_0$ and another  quintic $f_1$.   Let $(c_0, [f_0])\in \Gamma_{\mathbb L_1}$ be generic. Then
 
\par
(a)  \begin{equation}
{T_{c_0}\Gamma_{f_0}\over ker} \simeq H^0(c_0^\ast(T_{X_0})).\end{equation}
where $ker$ is a line in $T_{c_0}\Gamma_{f_0}$. 

\par
(b) 
\begin{equation} dim( T_{(c_0, [f_0])}\Gamma_{\mathbb L_1})=dim( T_{c_0}\Gamma_{f_0})+1, 
\end{equation}
and furthermore
\begin{equation} dim( T_{c_0}P(\Gamma_{\mathbb L_1}))=dim( T_{c_0}\Gamma_{f_0})+1, 
\end{equation}
\par

(c) If $dim( T_{c_0}P(\Gamma_{\mathbb L_1}))$=5, 
then\par

(1) $c_0$ is an immersion, \par
 (2) 
\begin{equation}
H^1(N_{c_0/X_0})=0.\end{equation}
\end{lemma}

\bigskip

\begin{proof}
(a). Let $a_i(c, f), i=0, \cdots, 5d$ be the coefficients of 
polynomial $f(c(t))$ in parameter $t$. Then the scheme
$$\Gamma$$ in $M\times \mathbf P^4$ is defined by homogeneous
polynomials $$a_i(c, f)=0, i=0, \cdots, 5d, \ locally\ around\ (c_0, [f_0]).$$ 
 Let $\alpha\in T_{c_0}M$. The equations
\begin{equation} {\partial  a_i(c_0, f_0)\over \partial \alpha}=0, \ all\ i \end{equation}  by the definition,  are necessary and sufficient conditions
for $\alpha$ lying in $$T_{(c_0, [f_0])}\Gamma_{f_0}.$$ 
On the other hand there is an evaluation map $e$:
\begin{equation}\begin{array}{ccc}
M\times \mathbf P^1 &\rightarrow & \mathbf P^4\\
(c, t) &\rightarrow & c(t)
\end{array}\end{equation}
The differential map
$e_\ast$ gives a morphism $e_m$:
\begin{equation}\begin{array}{ccc}
T_{c_0}M &\stackrel{e_m}\rightarrow & H^0(c_0^\ast( T_{\mathbf P^4}))\\
\alpha &\rightarrow & e_\ast (\alpha)
\end{array}\end{equation}

Suppose there is an $\alpha$ such that
$e_\ast (\alpha)=0$. We may assume $c_0$ is a map
$$\mathbb C^1\to \mathbb C^5-\{0\}.$$
 
Since $c_0$ is birational to its image, there is a Zariski open set  $$U_{\mathbf P^1}\subset \mathbb C^1\subset \mathbf P^1$$
and an open set $$V\subset \{(c_0(t))\}\subset \mathbb C^{5}-\{0\}$$
such that $c_0'|_{U_{\mathbf P^1}}$ is an isomorphism
$$U_{\mathbf P^1} \to V,$$ where
$c_0'$ is the map from $U_{\mathbf P^1}$ to $\mathbb C^{5}$ induced from $c_0$.
Due to the equation $e_\ast (\alpha)=0$, on $T_{t}U_{\mathbf P^1}$
$$(\alpha_0(t), \cdots, \alpha_4(t))=\lambda(t) c_0(t)$$
on $V$ (at each point $(c_0(t), \cdots, c_4(t))$ of $V$)
where $\lambda(t)$ lies in $\mathcal O(U_{\mathbf P^1})$. Because
$(\alpha_0(t), \cdots, \alpha_4(t))$ is parallel to
$c_0(t)$ at all points $t\in \mathbf P^1$, $\lambda$ can be extended to $\mathbf P^1$. Hence
$\lambda(t)$ is in $\mathcal O(\mathbf P^1)$. So it is a constant (independent of $t$).
Therefore $\alpha\in \mathbb C^{5d+5}$ is parallel to $$c_0\neq 0\in \mathbb C^{5d+5}.$$
This shows that 
$$dim(ker(e_m))=1.$$
(this does not hold for   a multiple cover map $c_0$ since the isomorphism $U_{\mathbf P^1} \to V$ does not exist for a multiple cover map.).  By the dimension count, $e_m$ must be surjective.
 
For any $\alpha\in c_0^\ast( T_{\mathbf P^4})$, 
$\alpha\in c_0^\ast(T_{X_0})$ if and only if
\begin{equation} {\partial  f_0(c_0(t))\over \partial \alpha}=0, \end{equation} 
for generic $t\in \mathbf P^1$. Notice equations (3.8) and (3.11) are exactly the same. Therefore
$e_m$ induces an isomorphism 
\begin{equation} \begin{array}{ccc} {T_{c_0}\Gamma_{f_0}\over ker(e_m)} &\stackrel{e_m}\rightarrow &  H^0(c_0^\ast(T_{X_0}))
\end{array} \end{equation}
This proves part (a).\par

(b). Let $t_1, \cdots, t_{5d+1}\in \mathbf P^1$ be $5d+1$ distinct points of $\mathbf P^1$.
Notice  $c^\ast(f)=0$ if and only if $c^\ast(f)|_{t_i}=0$ for all $i$. 
Then $$\Gamma_{\mathbb L_1}$$ is an irreducible component of 
\begin{equation} \{ (c, f)\in M\times \mathbb L_1: c^\ast(f)|_{t_i}=0, i=1, \cdots, 5d+1\}. \end{equation}
surjective to $T_{f_0}S$ in first order. 
Let $\alpha\in T_{c_0}M$ and $\beta=\overrightarrow f_1$. Let  $\epsilon\in \mathbb C$, then
 $\epsilon \beta\in T_{[f_0]}\mathbb L_1$. 
Then the Zariski tangent space $T_{c_0}\Gamma_{\mathbb L_1}$ is
$$\{(\alpha, \epsilon \beta)\in T_{(c_0, [f_0])}(M\times \mathbb L_1):
\epsilon f_1(c_0(t_i))+{\partial f_0(c_0(t_i))\over \partial \alpha} =0, i=1, \cdots, 5d+1\}$$
Because $[f_0]$ is a generic point of $\mathbb L_1$, the map (3.1) is surjective.  Thus
there is $\alpha_0\in T_{c_0}(M)$ such that 
$$ f_1(c_0(t))+{\partial f_0(c_0(t))\over \partial \alpha_0}=0$$ for all $t\in\mathbf P^1$. 
Thus $T_{(c_0, [f_0])}\Gamma_{\mathbb L_1}$ is isomorphic to
$$\{(\alpha, \epsilon)\in\mathbb C^{5d+6}:
{\partial f_0(c_0(t_i))\over \partial (\alpha-\epsilon\alpha_0)}=0, i=1,\cdots, 5d+1\}.$$

Notice that the subspace with $\epsilon=0$, 
$$\{\alpha\in T_{c_0}M:
{\partial f_0(c_0(t_i))\over \partial \alpha} =0, i=1,\cdots, 5d+1\},$$
is  the  tangent space of the scheme
$$\Gamma_{f_0}.$$
Thus the dimensions of them differ by $1$. 

\par

Next we prove the assertion for the scheme-theoretical image $P(\Gamma_{\mathbb L_1})$. 
Using the new expression of $\Gamma_{\mathbb L_1}$, $P(\Gamma_{\mathbb L_1})$ is defined
by polynomials
\begin{equation}
f_1(c(t_i))f_0(c(t_j))-f_1(c(t_j))f_0(c(t_i))=0
, 1\leq i, j\leq 5d+1.\end{equation}
Assume non of $t_i, i=1, \cdots, 5d+1$ is a zero of
$f_1(c_0(t))=0$. Then there exists an open set $U$ of $M$ around $c_0$ such that
\begin{equation} P(\Gamma_{\mathbb L_1})\cap U\end{equation}
is defined by
$5d$ equations
\begin{equation}
f_1(c(t_{5d+1}))f_0(c(t_j))-f_1(c(t_j))f_0(c(t_{5d+1}))=0, j=1, \cdots, 5d.\end{equation}
We'll denote $U$ by $M$.
Then the Zariski tangent space $T_{c_0} P(\Gamma_{\mathbb L_1})$ is defined by
\begin{equation}\begin{array} {cc} &
  f_1(c_0(t_{5d+1})){\partial f_0(c_0(t_j))\over \partial \alpha}-f_1(c_0(t_j)){\partial f_0(c_0(t_{5d+1}))\over \partial \alpha}=0, \\
& j=1, \cdots, 5d.\end{array}
\end{equation}

where $\alpha\in T_{c_0}M$.
This is the same as
\begin{equation}
{\partial f_0(c_0(t_j))\over \partial \alpha}-{f_1(c_0(t_j))\over f_1(c_0(t_{5d+1})) }{\partial f_0(c_0(t_{5d+1}))\over \partial \alpha}=0
\end{equation}

Next we view ${\partial f_0(c_0(t_j))\over \partial \alpha}$ as an element in
$$(T_{c_0}M)^\ast.$$
If $${\partial f_0(c_0(t_i))\over \partial \alpha}, i=1, \cdots, 5d+1$$ are linearly independent, then
$dim(T_{c_0}\Gamma_{f_0})=4$ and by (3.17), $$dim(T_{c_0}P(\Gamma_{\mathbb L_1}))=5. $$ The lemma is proved.
If $${\partial f_0(c_0(t_i))\over \partial \alpha}, i=1, \cdots, 5d+1$$ are linearly dependent, there are two cases:\par
(1) All solutions $\alpha_0$  to (3.17) satisfy $${\partial f_0(c_0(t_{5d+1}))\over \partial \alpha}=0.$$
(2) Some solutions $\alpha_0$ to (3.17) do not satisfy $${\partial f_0(c_0(t_{5d+1}))\over \partial \alpha}=0.$$
The case (1) is false. Because if $${\partial f_0(c_0(t_{5d+1}))\over \partial \alpha_0}=0$$
for all solutions to (3.17), 
then ${\partial f_0(c_0(t_j))\over \partial \alpha_0}=0$ for all $j=1, \cdots, 5d+1$. Hence all
solutions $\alpha_0$ to (3.17) must be  sections of $E$ (see definition 1.5 for $E$). As we know 
$$<\overrightarrow f_1>$$ is a solution to (3.17), but 
$${\partial f_0(c_0(t))\over \partial <\overrightarrow f_1>}=-f_1(c_0(t))$$ 
which means $<\overrightarrow f_1>$ is not a section of $E$. This is a contradiction. 

\par
In case (2), the solutions $\alpha_0$ to (3.17) must either satisfy 
$${\partial f_0(c_0(t_j))\over \partial \alpha}=0$$ for all $j=1, \cdots, 5d+1$, in which case, they are sections 
of $E$, or are uniquely expressed as
\begin{equation}
{\partial f_0(c_0(t_j))\over \partial \alpha_0}={f_1(c_0(t_j))\over f_1(c_0(t_{5d+1})) }{\partial f_0(c_0(t_{5d+1}))\over \partial \alpha_0}, 
j=1, \cdots, 5d+1
\end{equation}
which are exactly $<\overrightarrow f_1>$. 
Thus the vector $<\overrightarrow f_1>$ offers another dimension to 
$$dim(T_{c_0}P(\Gamma_{\mathbb L_1})).$$
So 
\begin{equation}
T_{c_0}P(\Gamma_{\mathbb L_1})=T_{c_0}\Gamma_{f_0}\oplus \mathbb C<\overrightarrow f_1>_M,
\end{equation}
where $<\overrightarrow f_1>_M\in e_m^{-1} (<\overrightarrow f_1>)$. 
The lemma is proved.
\bigskip

(c) If  $dim( T_{c_0}P(\Gamma_{\mathbb L_1}))$=5, then by part (a), (b), 
\begin{equation} 
dim(H^0(c_0^\ast(T_{X_0})))=3.
\end{equation}
Now we consider it from a different point of view.
Because $c_0$ is a birational map to its image,  there are finitely many points $t_i\in\mathbf P^1$ where the differential map

\begin{equation}\begin{array}{ccc}
(c_0)_\ast: T_{t_i}\mathbf P^1 &\rightarrow & T_{c_0(t_i)}\mathbf P^4
\end{array}\end{equation}
is not injective. Assume its vanishing order at $t_i$ is $m_i$ .  Let
\begin{equation}
m=\sum_i m_i.
\end{equation}
Let $s(t)\in H^0(\mathcal O_{\mathbf P^1}(m))$ such that $$div(s(t))=\Sigma_i m_it_i.$$

The sheaf morphism $(c_0)_\ast$ is injective and induces a composition morphism $\xi_s$ of sheaves
\begin{equation}\begin{array}{ccccc}
T_{\mathbf P^1} &\stackrel{(c_0)_\ast}  \rightarrow & c_0^\ast(T_{X_0})\otimes \mathcal I_{div(s)} & \stackrel{1\over s(t)}\rightarrow &
c_0^\ast(T_{X_0})\otimes \mathcal O_{\mathbf P^1}(-m),\end{array}
\end{equation}
where $\mathcal I_{div(s)}$ is the ideal sheaf of $div(s)$. 
It is easy to see that the induced bundle morphism $\xi_b$ is injective. Let 
\begin{equation}
N_m={c_0^\ast(T_{X_0})\otimes \mathcal O_{\mathbf P^1}(-m)\over \xi_b(T_{\mathbf P^1})}.
\end{equation}
Then
\begin{equation}
dim(H^0(N_m))=dim(H^0(c_0^\ast(T_{X_0})\otimes \mathcal O_{\mathbf P^1}(-m)))-3.
\end{equation}
On the other hand, three dimensional automorphism group of $\mathbf P^1$ gives a rise
to a 3-dimensional subspace  $Au$ of  $$H^0(c_0^\ast(T_{X_0})).$$
By (3.21), $Au=H^0(c_0^\ast(T_{X_0}))$. Over each point $t\in \mathbf P^1$, $Au$ spans 
a one dimensional subspace. Hence 
\begin{equation} c_0^\ast(T_{X_0})=\mathcal O_{\mathbf P^1}(2)\oplus \mathcal O_{\mathbf P^1}(-k_1)\oplus \mathcal O_{\mathbf P^1}(-k_2),
\end{equation}
where $k_1, k_2$ are some positive integers.
This implies that
\begin{equation} dim (H^0(c_0^\ast(T_{X_0})\otimes \mathcal O_{\mathbf P^1}(-m)))=dim( H^0(\mathcal O_{\mathbf P^1}(2-m)).
\end{equation}
Then
\begin{equation} dim (H^0(c_0^\ast(T_{X_0})\otimes \mathcal O_{\mathbf P^1}(-m)))=3-m.\end{equation}

Since $dim(H^0(N_m))\geq 0$, by (3.26), $-m\geq 0$. 
By the definition of $m$, $m=0$.  Hence $c_0$ is an immersion. 
\par
Next we prove (2). Notice that $(c_0)_\ast (T_{\mathbf P^1})$ is a subbundle generated by
global sections. It must be the $\mathbf O_{\mathbf P^1}(2)$ summand in (3.27) because $k_1, k_2$ are positive.
Therefore 
\begin{equation}
N_{c/X_0}\simeq \mathcal O_{\mathbf P^1}(-k_1)\oplus \mathcal O_{\mathbf P^1}(-k_2).
\end{equation}
Since $deg(c_0^\ast(T_{X_0}))=0$,  $k_1=k_2=1$.

Therefore 
\begin{equation}
H^1(N_{c_0/X_0})=0.\end{equation}

\end{proof}

\bigskip

Now we can describe the case when this is based on 2 dimensional plane in $S$. 
Recall $\mathbb L$ is an open set of a plane spanned by $f_0, f_1, f_2$, and $$\mathbb L_1\subset \mathbb L$$ is a pencil 
containing $[f_0]$.  
\bigskip

\begin{lemma}
For generic $c_g\in P(\Gamma_{\mathbb L_1})\subset P(\Gamma_{\mathbb L})$,
\begin{equation} 
dim( T_{c_g}P(\Gamma_{\mathbb L}))=dim( T_{c_g}P(\Gamma_{\mathbb L_1}))+1
\end{equation}

\end{lemma}

\bigskip

\begin{proof} 

Consider an open set $U_{P(\Gamma_{\mathbb L})}$ of ${P(\Gamma_{\mathbb L})}$ centered at $c_g$ such that
vectors in $\mathbb C^3$, 
$$\begin{array}{cc} & \biggl( f_2(c(t_1)), f_1(c(t_1)),  f_0(c(t_1))\biggr)\\&
\biggl(f_2(c(t_2)),  f_1(c(t_2)), f_0(c(t_2))\biggr)\end{array}$$
and vectors in $\mathbb C^{5d+1}$
\begin{equation} 
\left (\begin{array}{cc}  &f_2(c(t_{5d+1}) \\&
\vdots \\ &
f_2(c(t_2)) \\
 & f_2(c(t_1)) \end{array}\right),  \left (\begin{array}{cc}  &f_1(c(t_{5d+1}) \\&
\vdots \\ &
f_1(c(t_2)) \\
 & f_1(c(t_1)) \end{array}\right)
\end{equation}
are linearly independent for all $c\in U_{P(\Gamma_{\mathbb L})}$. 
If $$\left|  \begin{array}{ccc} f_2(c(t_1)) & f_1(c(t_1)) & f_0(c(t_1))\\
f_2(c(t_2)) & f_1(c(t_2)) & f_0(c(t_2))\\
f_2(c(t_i)) & f_1(c(t_i)) & f_0(c(t_i))\end{array}\right|=0,$$

$$\biggl(f_2(c(t_i)),  f_1(c(t_i)), f_0(c(t_i))\biggr)$$ for all $i=1, \cdots, 5d+1$ are linear combinations of 
$$\begin{array}{cc} & \biggl( f_2(c(t_1)), f_1(c(t_1)),  f_0(c(t_1))\biggr)\\& 
\biggl(f_2(c(t_2)),  f_1(c(t_2)), f_0(c(t_2))\biggr).\end{array}$$
Hence 
$$\left|  \begin{array}{ccc} f_2(c(t_i)) & f_1(c(t_i)) & f_0(c(t_i))\\
f_2(c(t_j)) & f_1(c(t_j)) & f_0(c(t_j))\\
f_2(c(t_l)) & f_1(c(t_l)) & f_0(c(t_l))\end{array}\right|=0,$$
for all $i, j, l$ between $1$ and $5d+1$. 

Since $$\left|  \begin{array}{ccc} f_2(c(t_i)) & f_1(c(t_i)) & f_0(c(t_i))\\
f_2(c(t_j)) & f_1(c(t_j)) & f_0(c(t_j))\\
f_2(c(t_l)) & f_1(c(t_l)) & f_0(c(t_l))\end{array}\right|=0,$$
for all $i, j, l$ between $1$ and $5d+1$ define $U_{P(\Gamma_{\mathbb L})}$, 
$5d-1$ equations
$$\left|  \begin{array}{ccc} f_2(c(t_1)) & f_1(c(t_1)) & f_0(c(t_1))\\
f_2(c(t_2)) & f_1(c(t_2)) & f_0(c(t_2))\\
f_2(c(t_i)) & f_1(c(t_i)) & f_0(c(t_i))\end{array}\right|=0,$$
for $i=3, \cdots, 5d+1$
define the scheme $U_{P(\Gamma_{\mathbb L})}$. 
Then by the definition there is $f_g\in \mathbb L$ such that
\begin{equation} 
c_g^\ast(f_g)=0. 
\end{equation}

We denote $U_{\mathbb L}$ by  $\mathbb L$.  Then
 $P(\Gamma_{\mathbb L}))$ is defined by the polynomial equations

\begin{equation} 
\left|  \begin{array}{ccc} f_2(c(t_i)) & f_1(c(t_i)) & f_0(c(t_i))\\
f_2(c(t_1)) & f_1(c(t_1)) & f_0(c(t_1))\\
f_2(c(t_2)) & f_1(c(t_2)) & f_0(c(t_2))\end{array}\right| =0\end{equation}
 $i=3, \cdots, 5d+1$.
\par

We may assume that $c_g=c_0$ lies in $f_0$ (by choosing appropriate basis $f_0, f_1, f_2$ of $\mathbb L$).  By  (3.35),
the Zariski tangent space $T_{c_0}P(\Gamma_{\mathbb L})$ is defined by equations
\begin{equation} 
\left|  \begin{array}{ccc} f_2(c_0(t_i)) & f_1(c_0(t_i)) & {\partial f_0(c_0(t_i))\over\partial \alpha}\\
f_2(c_0(t_1)) & f_1(c_0(t_1)) & {\partial f_0(c_0(t_1))\over\partial \alpha}\\
f_2(c_0(t_2)) & f_1(c_0(t_2)) & {\partial f_0(c_0(t_2))\over\partial \alpha}\end{array}\right| =0\end{equation}
 where $i=3, \cdots, 5d+1$, $\alpha\in T_{c_0}P(\Gamma_{\mathbb L})$.
This is equivalent to that 
the column vectors of 
\begin{equation} 
\left (\begin{array}{ccc} f_2(c_0(t_{5d+1}) & f_1(c_0(t_{5d+1})) & {\partial f_0(c_0(t_{5d+1}))\over\partial \alpha_0}\\
\vdots & \vdots &\vdots\\
f_2(c_0(t_2)) & f_1(c_0(t_2)) & {\partial f_0(c_0(t_{2}))\over\partial \alpha_0}\\
f_2(c_0(t_1)) & f_1(c_0(t_1)) & {\partial f_0(c_0(t_{1}))\over\partial \alpha_0}\end{array}\right) \end{equation}
 are linearly dependent,  for some  $\alpha_0\in T_{c_0}P(\Gamma_{\mathbb L})$. 
Since 
\begin{equation} 
\left (\begin{array}{cc}  &f_2(c_0(t_{5d+1}) \\&
\vdots \\ &
f_2(c_0(t_2)) \\
 & f_2(c_0(t_1)) \end{array}\right),  \left (\begin{array}{cc}  &f_1(c_0(t_{5d+1}) \\&
\vdots \\ &
f_1(c_0(t_2)) \\
 & f_1(c_0(t_1)) \end{array}\right)
\end{equation}
are linearly independent,
there are complex numbers $\epsilon_i, i=1, 2$ such that
\begin{equation} 
\left (  \begin{array}{c}  {\partial f_0(c_0(t_{5d+1}))\over\partial \alpha_0}\\
 \vdots\\
  {\partial f_0(c_0(t_{2}))\over\partial \alpha_0}\\
 {\partial f_0(c_0(t_{1}))\over\partial \alpha_0}\end{array}\right)=\epsilon_2 \left (  \begin{array}{c}  f_2(c_0(t_{5d+1})\\
 \vdots\\
 f_2(c_0(t_2))\\
 f_2(c_0(t_1))\end{array}\right)+\epsilon_1 \left (  \begin{array}{c} f_1(c_0(t_{5d+1})\\
 \vdots\\
f_1(c_0(t_2))\\
 f_1(c_0(t_1))\end{array}\right)
\end{equation}

Let $f_3=\sum_{i=1}^2 \epsilon_i f_i$ where $\epsilon_i$ are fixed.   We may assume that $\mathbb L_1$ is the pencil containing $f_0, f_3$.  The equation (3.39) becomes
$$ {\partial f_0(c_0(t))\over \partial \alpha_0}-f_3(c_0(t))=0, for\ all\ t\in\mathbf P^1.$$

Then just as in the (3.20), 
\begin{equation}
T_{c_0}\Gamma_{\mathbb L}\simeq T_{c_0}\Gamma_{f_0}\oplus \mathbb C<\overrightarrow{f_1}>_M\oplus \mathbb C<\overrightarrow{f_2}>_M.
\end{equation}
and \begin{equation}
T_{c_0}\Gamma_{\mathbb L_1}\simeq T_{c_0}\Gamma_{f_0}\oplus \mathbb C<\overrightarrow{f_3}>_M.
\end{equation}

Then the lemma follows.

\end{proof}

\bigskip

\section{Space of rational curves, $M$}

\bigskip
The basis of this paper is the linear model of stable moduli, which begins with the projectivization $\mathbf P(M)$. The space $M$ is an affine space
$\mathbb C^{5d+5}$, therefore is very simple.  But we are interested in some subschemes which are not trivial at all. Our idea is to introduce various
analytic coordinates of each copy $\mathbb C^{d+1}$ in $\mathbb C^{5d+5}$. The purpose of these coordinates is to reveal the local higher-order structures of
the subschemes.  Those structures will be broken-up to multiple pieces of manageable sizes ( in terms of Jacobian matrices), and there are no unified coordinates even at a fixed point to handle all pieces. Thus in this section we introduce a couple of coordinates systems that will be used interchangeably.   But  these coordinates are not for the $M$, rather for the 
successive blow-ups of $M$.  Analytic neighborhoods of these blow-ups  can be identified with neighborhoods $\mathbb C^{5d+5}$. 
So for now we just borrow $M$ to set-up these coordinates for the later calculations. 
\bigskip

However readers may skip this section because without section 5 technical preparation here may seem to be aimless. \bigskip

Let $\tilde c_2=(\tilde c_2^0, \cdots, \tilde c_2^4)\in M$ with $$\tilde c_2^i\in H^0(\mathcal O_{\mathbf P^1}(d))-\{0\}, i=0, \cdots, 4.$$
We may assume $t\in \mathbb C\subset \mathbf P^1$. 
We assume $\tilde c_2^i(t)=0, i=0, \cdots, 4$ have  $5d$ distinct  zeros $$\tilde \theta_i^j, for\ i\leq 4.$$
Then each component, $ H^0(\mathcal O_{\mathbf P^1}(d))$ of 
$$ M=H^0(\mathcal O_{\mathbf P^1}(d))^{\oplus 5}$$
 has local analytic coordinates \begin{equation}
r_i, \theta_i^j, j=1, \cdots, d, for \ r_i\neq 0\end{equation}
 (for each $i=0, 1, 2, 3, 4$) around
$\tilde c_2^i$ such that
\begin{equation} c^i(t)=r_i \prod _{j=1}^d (t-\theta_i^j).  
\end{equation}
 
Let  coordinates values for $\tilde c_2$ be
$$r_l=y_l, \theta_i^j=\tilde\theta_i^j,  l=0, \cdots, 4,  i=0, \cdots, 4,  j=1, \cdots, d.$$
Let $q$ be a generic, homogeneous quadratic polynomial in $z_0, \cdots, z_4$.

Let 
\begin{equation} 
h(c, t)=\delta_1q(c(t))+\delta_2 c_3(t) c_4(t).
\end{equation}
for $c\in M$, where $\delta_i, i=1, 2$ are two none zero complex numbers.
   Let $\beta_1, \cdots, \beta_{2d}$ be the zeros of 
$h(\tilde c_2, t)=0$. Also let  
$$h(c, t)=\xi \prod _{i=1}^{2d}(t-\epsilon_i).$$
It is clear that $$\begin{array}{c}
\xi=\delta_1 q(r_0, r_1, r_2, r_3, r_4)+\delta_2 r_3r_4, \ and\\
\epsilon_i\ are \ analytic\ functions\ of\ c.
\end{array}$$
Let the corresponding value of $\xi$ at $\tilde c_2$ be $\xi^0$. 
By the genericity of $q$, we may assume $\beta_i, i=1, \cdots, 2d$ are distinct and non-zeros . 
Furthermore we assume $\beta_i$ are distinct for $q=z_1z_2$ and generic $\delta_i$.

\bigskip

\begin{proposition}
Let $U_{\tilde c_2}\subset M$ be an analytic neighborhood of $\tilde c_2$. 

\par
(a) Let \begin{equation}\begin{array}{ccc}
\varrho: U_{\tilde c_2} &\rightarrow &  \mathbb C^{5d+5}
\end{array}\end{equation}
be a regular map that is defined
by 
\begin{equation}\begin{array}{cc} &
 (\theta_0^1, \cdots, \theta_4^d, r_0, r_1, r_2, r_3, r_4) \\
 &\xdownarrow{0.3cm}\scriptstyle{\varrho} \\ &
(\theta_0^1, \cdots, \theta_2^d, \epsilon_1, \cdots, \epsilon_{2d}, 
r_0, \cdots, r_3, \xi).
\end{array}\end{equation} 

Then $\varrho$ is an isomorphism to its image.

\par
(b) Let \begin{equation}\begin{array}{ccc}
\varrho': U_{\tilde c_2} &\rightarrow &  \mathbb C^{5d+5}
\end{array}\end{equation}
be a regular map that is defined
by 
\begin{equation}\begin{array}{cc} &
(\theta_0^1, \cdots, \theta_4^d, r_0, r_1, r_2, r_3, r_4) \\
 &\xdownarrow{0.3cm}\scriptstyle{\varrho'}\\ &
(\theta_0^1, \cdots, \theta_2^d, \epsilon_1, \cdots, \epsilon_{2d}, 
r_0, \cdots, r_3, r_4).
\end{array}\end{equation} 
Then $\varrho'$ is an isomorphism to its image.
\end{proposition}
\bigskip

\begin{proof}
It suffices to prove the differential of $g$ at $\tilde c_2$ is an isomorphism for a SPECIFIC $q$. So we assume  that 
$$ \delta_1=\delta_2=1, q=z_1z_2$$
This is a straightforward calculation of the Jacobian of $g$. We may still assume that $\beta_i, i=1, \cdots, 2d$ are distinct. Using
the composition of two isomorphisms, we obtain that 
the Jacobian 
\begin{equation}
{\partial g(\tilde \theta_0^1, \cdots, \tilde \theta_2^d, y_0,\cdots, y_3, \xi^0,   \beta_1,\cdots, \beta_{2d})
\over \partial (\theta_0^1, \cdots, \theta_2^d, r_0, r_1, r_2, r_3, r_4, \theta_3^1, \cdots, \theta_4^d)}
\end{equation} is equal to 

\begin{equation}
a\cdot {\partial \xi\over \partial r_4}|_{\tilde c_2} \cdot Ja
\end{equation}
where $a$ is some non-zero number, ${\partial \xi\over \partial r_4}|_{\tilde c_2}$ is also non-zero and $Ja$ is another Jacobian

\begin{equation}
Ja=\left|\begin{array}{ccc}
{\partial  h(c, \beta_1)\over \partial \theta_3^1}  &\cdots  & {\partial h(c, \beta_{1})\over \partial\theta_4^d }\\
\vdots & \vdots & \vdots\\
{\partial  h(c, \beta_{2d})\over \partial \theta_3^1}  &\cdots  & {\partial h(c, \beta_{2d})\over \partial\theta_4^d }
\end{array}\right|_{\tilde c_2}
\end{equation}
Let $T_i, i=0, d$   be the determinant
\begin{equation}
\left|\begin{array}{ccc}
\beta_{i+1}  & \cdots & \beta_{i+1}^d\\
\vdots & \vdots & \vdots \\
\beta_{i+d}  & \cdots & \beta_{i+d}^d
\end{array}\right|.
\end{equation}
Then we compute the determinant to have 
\begin{equation}
Ja=(-1)^d\delta_2 T_0T_d \prod_{i=1}^d (\tilde c_2^3(\beta_{d+i})\tilde c_2^4(\beta_{i}) -\tilde c_2^3(\beta_{i}) \tilde c_2^4(\beta_{d+i})).\end{equation}
Since $\beta_i$ are distinct and non-zeros, $$T_0\neq 0, T_d\neq 0.$$
Since ${\tilde c_2^3(t)\over \tilde c_2^4(t)}$ is a rational function and
$$deg(\tilde c_2^3(t))=deg(\tilde c_2^4(t))=d$$
then we can always arrange the index of $\beta_i$ such that each number 
\begin{equation} ({\tilde c_2^3(\beta_{d+i})\over \tilde c_2^4(\beta_{d+i})}-{\tilde c_2^3(\beta_{i})\over \tilde c_2^4(\beta_{i})})
\end{equation} is not zero. Hence 
$$\prod_{i=1}^d (\tilde c_2^3(\beta_{d+i})\tilde c_2^4(\beta_{i}) -\tilde c_2^3(\beta_{i}) \tilde c_2^4(\beta_{d+i}))\neq 0.$$

Thus $Ja$ is non-zero.   Therefore 
\begin{equation}
{\partial g(\tilde \theta_0^1, \cdots, \tilde \theta_2^d, y_0,\cdots, y_3, \xi^0,  h(\tilde c_2, \beta_1),\cdots, h(\tilde c_2,  \beta_{2d}))
\over \partial (\theta_0^1, \cdots, \theta_2^d, r_0, r_1, r_2, r_3, \xi, \theta_3^1, \cdots, \theta_4^d)}\neq 0
\end{equation}
 The proof of part (b) is the same as for part (a). 
We complete the proof. 

\end{proof}

\bigskip

\begin{definition}
By proposition 4.1,  both
\begin{equation} \theta_0^1, \cdots, \theta_2^d, r_0, \cdots, r_4, \epsilon_1, \cdots, \epsilon_{2d}
\end{equation} and 
\begin{equation} \theta_0^1, \cdots, \theta_2^d, r_0, \cdots, r_3, \xi, \epsilon_1, \cdots, \epsilon_{2d}
\end{equation}
are the local analytic coordinates of $M$ around $\tilde c_2$, and
$\tilde c_2$ corresponds to the coordinate values
\begin{equation}\begin{array} {c}
 \theta_i^j=\tilde \theta_i^j, i=0, 1,2, j=1, \cdots, d\\
r_l=y_l\neq 0, l=0, \cdots, 4\\
\epsilon_i=\beta_i, i=1, \cdots, 2d
\end{array}\end{equation}
and $\xi^0$. 

\end{definition}

\bigskip

\section{Differential sheaf }
\bigskip

In this section, we prove theorem 1.1, i.e. \bigskip

\begin{equation} H^1(N_{c_0/X_0})=0\end{equation}
at generic $(c_0, [f_0])\in \Gamma$.

\bigskip

\bigskip

The following lemma is also a local expression for the calculation in the section 5.1 below. 
Choose  homogeneous coordinates  $[z_0, \cdots, z_4]$ for $\mathbf P^4$. 
Let 
\begin{equation} 
f_3=z_0z_1z_2(\delta_1q+ \delta_2z_3 z_4).
\end{equation}
where $\delta_i$ are two non-zero complex numbers, and $q$ is a generic, quadratic homogeneous polynomial   in $z_0, \cdots, z_4$. 
Let 
$\tilde c_2\in M$ and 
$$f_3(\tilde c_2( t))\neq 0.$$
We denote the zeros of $\tilde c_2^i(t)=0$ by $\tilde\theta_i^j$ and zeros of
\begin{equation}
(\delta_1q+ \delta_2z_3 z_4|_{\tilde c_2(t)})=0
\end{equation}
by $\beta_i, i=1, \cdots, 2d$. 
We assume $\tilde \theta_i^j$ are distinct. 

\bigskip

\begin{lemma} 
Let $t_1, \cdots, t_{5d}$ be the zeros of  $f_3(\tilde c_2( t))$, i.e. they are equal to
$\tilde\theta_i^j, i\leq 2$ and $\beta_i$.    Recall 
in definition 4.2, 
$$ \theta_0^1, \cdots, \theta_2^d, r_0, \cdots, r_3, \xi,  \epsilon_1, \cdots, \epsilon_{2d}
$$ are analytic coordinates of $M$ around the point $\tilde c_2$.

Then \par

(a) the Jacobian matrix

\begin{equation}\begin{array} {c} J( \tilde c_2)\\
\|\\
\left(\begin{array}{cccccc}
{\partial f_3(\tilde c_2(t_{1}))\over \partial \theta_0^1}  &\cdots &
{\partial f_3(\tilde c_2(t_{1}))\over \partial \theta_2^d} & {\partial f_3(\tilde c_2(t_{1}))\over \partial \epsilon_1} &\cdots &
{\partial f_3(\tilde c_2(t_{1}))\over \partial \epsilon_{2d}}\\
{\partial f_3(\tilde c_2(t_{2}))\over \partial \theta_0^1}  &\cdots &
{\partial f_3(\tilde c_2(t_{2}))\over \partial \theta_2^d} & {\partial f_3(\tilde c_2(t_{2}))\over \partial \epsilon_1} &\cdots &
{\partial f_3(\tilde c_2(t_{2}))\over \partial \epsilon_{2d}}\\
{\partial f_3(\tilde c_2(t_{3}))\over \partial \theta_0^1}  &\cdots &
{\partial f_3(\tilde c_2(t_{3}))\over \partial \theta_2^d} & {\partial f_3(\tilde c_2(t_{3}))\over \partial \epsilon_1} &\cdots &
{\partial f_3(\tilde c_2(t_{3}))\over \partial \epsilon_{2d}}\\
\vdots&\vdots &\vdots &\vdots &\vdots&\vdots\\
\vdots&\vdots &\vdots &\vdots &\vdots&\vdots\\
{\partial f_3(\tilde c_2(t_{5d}))\over \partial \theta_0^1}  &\cdots &
{\partial f_3(\tilde c_2(t_{5d}))\over \partial \theta_2^d} & {\partial f_3(\tilde c_2(t_{5d}))\over \partial \epsilon_1} &\cdots &
{\partial f_3(\tilde c_2(t_{5d}))\over \partial \epsilon_{2d}}
\end{array}\right)\end{array}\end{equation}
is equal to 
a diagonal matrix $D$ whose diagonal entries are
\begin{equation}
{\partial f_3(\tilde c_2(t_1))\over \partial \theta_0^1}, \cdots, {\partial f_3(\tilde c_2(t_{3d}))\over \partial \theta_2^d}, 
{\partial f_3(\tilde c_2( t_{3d+1}))\over \partial \epsilon_1}, \cdots, 
{\partial f_3(\tilde c_2( t_{5d}))\over \partial \epsilon_{2d}}
\end{equation}
which are all non-zeros.\par
(b) For $i=1, \cdots, 5d$, $j=0, 1, 2, 3$, 
$${\partial f_3(\tilde c_2(t_i))\over \partial r_j}={\partial f_3(\tilde c_2(t_i))\over \partial \xi}=0.$$

\end{lemma}

\bigskip

\begin{proof}  Note $\tilde \theta_i^j$ are distinct and $\beta_i$ are also distinct by the genericity of $q$.  Thus the coordinates
in definition 4.2 exist. 
It suffices to show all non diagonal entries of (5.4) are zeros. We can rewrite 
\begin{equation}
f_3(c(t))=y \prod_{j=1}^{5d}(t-\alpha_j)
\end{equation}
around $\tilde c_2$.
Then $y, \alpha_j$ are all functions of the analytic coordinates in definition 4.2, 
\begin{equation}
\theta_i^j, \epsilon_1, \cdots, \epsilon_{2d}, y_l.
\end{equation}
More specifically $\theta_i^j, i\leq 2$ are exactly the $3d$ roots of $f_3(c(t))$ and $\epsilon_i$ are 
the parameters in the local ring, which are defined by the evaluations of $h(c, t)$ (a factor of $f_3(c(t))$) at its $2d$ roots. 
This directly implies that the Jacobian matrix  
\begin{equation}
{\partial (\alpha_1, \cdots,  \alpha_{5d})\over \partial (\theta_0^1, \cdots, \theta_2^d, \epsilon_1,\cdots,\epsilon_{2d})}
\end{equation}
is a  non-zero diagonal matrix when evaluated at $\tilde c_2$.  Since it is clear that
$$
{\partial f_3 (\tilde c_2(t_i))\over \partial \alpha_j}=0, i\neq j, 
$$
the non-diagonals of $J(\tilde c_2)$ are also $0$. 
Up to a non-zero constant multiple, the diagonal entries of $J(\tilde c_2)$ are 
$$\prod_{j\neq 1}(t_1-t_j), \prod_{j\neq 2}(t_2-t_j), \cdots, \prod_{j\neq 5d}(t_{5d}-t_j).$$
For part (b), we write down an expression of $f_3 (c(t))$, 
$$f_3 (c(t))=r_0 r_1r_2 \xi \prod_{i=0, j=1, l=1}^{i=2, j=d, l=2d }(t-\theta_i^j)(t-\epsilon_l).$$
Part (b) follows from this expression.
We complete the proof. 
\end{proof}

\bigskip

\subsection{Non-vanishing $5d-1$-form $\omega(M, \mathbf t)$ }
\quad\smallskip

The section 4 and lemma 5.1 are just the preparation for the proof. Section 3 is the one that is meant to dig into  this problem, but it is
only to the first order.  With only the first order results in section 3, we can't go far because Clemens' conjecture, we believe, touches upon  all orders of deformations of pairs.
The following lemma is a reflection of this philosophy.\bigskip

\begin{lemma} 

The $5d$-$1$ form $\omega(M, \mathbf t)$ defined in (1.8) is a non-zero form when it is evaluated at generic points of $P(\Gamma_{\Bbb L})$, i.e.
the reduction $\bar \omega(M, \mathbf t)$ in the module, 
$$H^0(\Omega_M\otimes \mathcal O_{P(\Gamma_{\Bbb L})})$$ is non zero.

\end{lemma}
\bigskip

This lemma is  proposition 1.3 in the introduction.
\bigskip

It suffices to prove lemma 5.2 for a special choice of $f_0, f_1, f_2$, and we only need to produce one point (any point) on $P(\Gamma_{\mathbb L})$, 
at which $\omega(M, \mathbf t)$ is non-vanishing.\par
So let $z_0, z_1, \cdots, z_4$ be  general homogeneous coordinates of
$\mathbf P^4$. Let $$f_2=z_0z_1z_2z_3z_4.$$  
Let 
$$f_1=z_0z_1z_2 q,$$
where $q$ is a generic quadratic homogeneous polynomial in $z_0, \cdots, z_4$. 
   Choose another generic $f_0$.  
We obtain an open set $\mathbb L_1$ of pencil through $f_0, f_2$, and an
 open set $\mathbb L$ of 2-dimensional plane containing $f_0, f_1, f_2$ in $S$.
 Let $P(\Gamma_{\mathbb L})$ and $P(\Gamma_{\mathbb L_1})$ be as defined in lemmas 3.2, 3.3.  We choose $P(\Gamma_{\mathbb L_1})$ to be irreducible, and to 
be contained in $P(\Gamma_{\mathbb L})$ for all generic  $q$ (simultaneously). 
We may assume a generic point $c=(c^0, \cdots, c^4)\in  P(\Gamma_{\mathbb L_1})$ does not have multiple zeros with coordinates planes, i.e.
$c^i(t)=0, i=0, \cdots, 4$ have $5d$ distinct roots.  This is because $c$ is a birational map to its image.   \par

Let $$c_2\in P(\Gamma_{\mathbb L_1})\subset P(\Gamma_{\mathbb L})$$
be a special point 
such that $f_2(c_2(t))=0$.    By the genericity of $f_0$ and $z_0$, 
we assume 
$$P(\Gamma_{\mathbb L_1})\not\subset \{c^0=0\}.$$
For the special point $c_2\in M$, there are two factors that would result in $\omega(M, \mathbf t)=0$:
(1) vanishing components $c^i_2$ where $c_2=[c_2^0, c^1_2, \cdots, c_2^4]$, \par
(2) common zeros of $\tilde c_2(t)^i=0$ for the blow-ups of $c_2$.\par
The second factor will become clearer after the first factor is resolved. Thus 
 let's start with the first case:
we assume $c_2=[0, c^1_2, \cdots, c_2^4]$ where
$c_2^i, i\neq 0$ are non-zero sections of $H^0(\mathcal O_{\mathbf P^1}(d))$.\footnote{ $c_2$ occurs as a ``boundary point" in the ``compactification" of the space of rational maps in our setting. It could be as ``bad" as a constant map later in the proof. So in general we should not regard $c_2$ as a map.} 
It is not difficult to see $\omega(M, \mathbf t)$ is zero at $c_2$ for the choice
of $\mathbb L$.    But we would like to show that $\omega(M, \mathbf t)$ is  not identically zero on  $P(\Gamma_{\mathbb L})$.    The technique is ``blow-up".  \par
Let's  first describe it in the most general term: we would
construct a birational map, the composition of successive blow-ups, 
\begin{equation}\begin{array}{ccc}
\pi: Y &\rightarrow &P(\Gamma_{\mathbb L})
\end{array}\end{equation}
which is an isomorphism on
$Y-\pi^{-1}(B')$ where $B'$ is a proper closed sub-scheme of $P(\Gamma_{\mathbb L})$ that contains $c_2$, and
$\pi^{-1}(B')$ is also a proper,  closed sub-scheme.
We then compute to find out that 
\begin{equation}
\pi^\ast(\omega(M, \mathbf t))= g \omega(M, \mathbf t)'
\end{equation}
on $Y$ where $g$ is a non-zero rational function on $Y$, and  $\omega(M, \mathbf t)'$ does not vanish at a point $p\in \pi^{-1}(B')$. Because $\pi$ is birational,
$\omega(M, \mathbf t)$ is a non-zero form restricted to  $P(\Gamma_{\mathbb L})$. The key to this assertion is that
$\omega(M, \mathbf t)'$ is non-zero at one point $p$. 
The process of blow-ups is very lengthy. In order to organize them,
we  divide them into two different types dealing with each factor mentioned above: \par
(1) blow-ups used to resolve the vanishing of sections $c_2^i, i=0, \cdots, 4$ in the rational map $c_2$. They are algebraic and  will be only used
in step 1 below.\par
(2) blow-ups used to resolve the multiple zeros of $\tilde c_2$ with coordinates planes. This
is the case where $c_2$ could be a constant map. They are analytic (however they can be replaced by algebraic blow-ups).  These will be  used in step 2 below. 

\par

Let's see the details.  
Let 
\begin{equation}
B_1\subset \mathbb C^{5d+5}
\end{equation}
be the subvariety 
that is equal to 
\begin{equation}
\{0\}\oplus  H^0(\mathcal O_{\mathbf P^1}(d))\oplus  H^0(\mathcal O_{\mathbf P^1}(d))\oplus  H^0(\mathcal O_{\mathbf P^1}(d))\oplus  H^0(\mathcal O_{\mathbf P^1}(d)).
\end{equation}
Let   \begin{equation}\begin{array}{ccc}
 \tilde \mathbb C^{5d+5} &\stackrel{\pi_1} \rightarrow & 
\mathbb C^{5d+5}\end{array}\end{equation}
 be the blow-up of $\mathbb C^{5d+5}$ along $B_1$. It is clear that
\begin{equation}\begin{array}{cc}
 & \tilde \mathbb C^{5d+5}\\
& \|\\ &
\tilde \mathbb C^{d+1}\times H^0(\mathcal O_{\mathbf P^1}(d))\times H^0(\mathcal O_{\mathbf P^1}(d))\times H^0(\mathcal O_{\mathbf P^1}(d))\times H^0(\mathcal O_{\mathbf P^1}(d))\end{array}
\end{equation}
where  $\tilde \mathbb C^{d+1}$ is the blow-up of $\mathbb C^{d+1}$ at the origin.

Now let $\tilde M$ and $\tilde P(\Gamma_{\mathbb L})$ ($\tilde P(\Gamma_{\mathbb L_1})$ ) be the strict transforms of $M$ and 
$P(\Gamma_{\mathbb L})$ ($P(\Gamma_{\mathbb L_1})$ )  respectively. Let
$ E_1$, $E_{1\mathbb L}$ be their exceptional divisors.
 Let $\tilde c_2\in \tilde P(\Gamma_{\mathbb L_1})$ be an inverse of $c_2$ under the map
\begin{equation}\begin{array}{ccc}
 \tilde P(\Gamma_{\mathbb L_1}) &\stackrel{\pi_1|_{ \tilde P(\Gamma_{\mathbb L_1})}} \rightarrow & 
P(\Gamma_{\mathbb L_1}).\end{array} 
\end{equation}
(such $\tilde c_2$ is independent of choice of $q$).  Let $U_{ \tilde \mathbb C^{5d+5}}$ be a neighborhood of $\tilde c_2$.
\bigskip

\begin{lemma} Let $(t_1, \cdots, t_{5d+1})\in Sym^{5d+1}(\mathbf P^1)$ be generic. Let
\begin{equation}\begin{array}{cc}  \psi_i=\left|  \begin{array}{cc} f_0(c(t_1)) &  f_2(c(t_1))\\
f_0(c(t_2)) & f_2(c(t_2))\end{array}\right| df_1(c(t_i))+\left|  \begin{array}{cc} f_2(c(t_1)) &  f_1(c(t_1))\\
f_2(c(t_2)) & f_1(c(t_2))\end{array}\right| df_0(c(t_i))& \\  +\left|  \begin{array}{cc} f_1(c(t_1)) &  f_0(c(t_1))\\
f_1(c(t_2)) & f_0(c(t_2))\end{array}\right| df_2(c(t_i))
\end{array}\end{equation}
be 1-forms on $M$ for $i=3, 4, \cdots, 5d+1$. Then
vectors
\begin{equation}\begin{array}{cc} &  \pi_1^\ast (\psi_i), i=3, \cdots, 5d+1\\
& \pi_1^\ast (df_0(c(t_1))), \pi_1^\ast (df_0(c(t_2))),\\
& \pi_1^\ast (df_1(c(t_1))), \pi_1^\ast (df_1(c(t_2))),\\
& \pi_1^\ast (df_2(c(t_1))), \pi_1^\ast (df_2(c(t_2)))
\end{array}\end{equation}
are linearly independent in the vector space $(T_{\tilde c}{\tilde{M}})^\ast$, i.e. they form a basis, where 
$\tilde c$ lies in a non-empty open set of  $ \tilde P(\Gamma_{\mathbb L})$, and $\tilde c\neq \tilde c_2$.  \footnote{  The vectors in lemma 5.3 come from the expansion of $\phi_i$: $$\begin{array}{cc}
  \phi_i= \left|  \begin{array}{cc} f_0(c(t_1)) &  f_2(c(t_1))\\
f_0(c(t_2)) & f_2(c(t_2))\end{array}\right| df_1(c(t_i))+\left|  \begin{array}{cc} f_2(c(t_1)) &  f_1(c(t_1))\\
f_2(c(t_2)) & f_1(c(t_2))\end{array}\right| df_0(c(t_i))& \\  +\left|  \begin{array}{cc} f_1(c(t_1)) &  f_0(c(t_1))\\
f_1(c(t_2)) & f_0(c(t_2))\end{array}\right| df_2(c(t_i))
+ \sum_{l=0, j=1}^{l=2, j=2} h_{lj}^i(c) df_l(c(t_j)) & \end{array}$$
where $h_{lj}^i$ are polynomials in $c$.}   

\end{lemma}
\bigskip

\begin{proof} of lemma 5.3:  The proof  is long. Thus we divide it into two steps.
It suffices to prove it for special $t_1, \cdots, t_{5d+1}$.  First let
$\tilde c_2$ be decomposed (according to (5.14)) to 
$$\tilde c_2=( \tilde c_2^0(t), \cdots, \tilde c_2^4(t))$$
where $$\tilde c_2^0(t)\in \mathbf P(H^0(\mathcal O_{\mathbf P^1}(d)), \tilde c_2^i(t))\in H^0(\mathcal O_{\mathbf P^1}(d)), i\neq 0.$$
and the zeros of all $\tilde c_2^i, i=0, \cdots, 4$ are $\tilde\theta_i^j$.  
\bigskip

{\bf Step 1}: 
 Suppose that
 \begin{equation} \tilde\theta_i^j, , i=0, \cdots, 4, j=1, \cdots, d,
\end{equation} are distinct. 
In this case, we 
let  $t_1, t_2$ be two points on $\mathbf P^1$ satisfying
\begin{equation}
\left|\begin{array}{cc} q|_{\tilde c_2(t_1)}  &  \tilde c_2^3(t_1) \tilde c_2^4 (t_1)\\
q|_{\tilde c_2(t_2)}& \tilde c_2^3(t_2) \tilde c_2^4 (t_2)\end{array}\right|=0.
\end{equation}
Let 
\begin{equation}
f_3=z_0z_1z_2( \delta_1 q+\delta_2 z_3z_4)
\end{equation}
where $$\delta_1=\left|\begin{array}{cc} f_0 (c_2(t_1))  &   f_2(\tilde c_2(t_1))\\
f_0 (c_2(t_2))& f_2(\tilde c_2(t_2))\end{array}\right|, \quad \delta_2=\left|\begin{array}{cc} f_1 (\tilde c_2(t_1))  &   f_0(c_2(t_1))\\
f_1 (\tilde c_2(t_2))& f_0( c_2(t_2))\end{array}\right|.$$
Then let 
$t_3, \cdots, t_{5d}$ be zeros of
\begin{equation}\begin{array}{c} \bigl ({1\over r_0}\pi_1^{\ast} (f_3(z))\bigr)|_{\tilde c_2(t)}\\
\|\\
\tilde c_2^0(t) \tilde c_2^1(t)\tilde c_2^2(t) \biggl ( 
\delta_1
q|_{\tilde c_2(t)}+  \delta_2 \tilde c_2^{3}(t)\tilde c_2^4(t)\biggr)\\
\|\\
0.\end{array}
\end{equation}
other than $\tilde\theta_0^1, \tilde\theta_1^1$.  Let $t_{5d+1}$ be generic.  

Because $f_0,  q$ are generic,  
the zeros of 
\begin{equation}  \delta_1
q|_{\tilde c_2(t)}+  \delta_2 \tilde c_2^3(t)\tilde c_2^4(t)
\end{equation}
are distinct and non-zeros. 
Thus \begin{equation}
t_1, t_2, t_3, \cdots, t_{5d+1}
\end{equation}
are distinct.  We are going to use the coordinates in section 4. We start with the usual coordinates of $\tilde \mathbb C^{5d+5}$.
Let $c_i^j, i=0, \cdots, 4, j=0, \cdots, d$ be the coefficients of
five tuples of sections of $H^0(\mathcal O_{\mathbf P^1}(d))$. They are the affine coordinates of
$$\mathbb C^{5d+5}.$$
Each section of $\mathcal O_{\mathbf P^1}(d)$ in an analytic neighborhood excluding those with multiple zeros can be
written as 
\begin{equation}
c_i(t)=r_i\Pi_{j=1}^d (t-\theta_i^j).
\end{equation}
for $i\leq 4$. 
Then $r_m, \theta_i^j$ are local analytic coordinates for an analytic open set $U_{\tilde \mathbb C^{5d+5}}$ of the blow-up 
\begin{equation}  \tilde \mathbb C^{5d+5}
\end{equation} 
centered around $\tilde c_2\in \tilde \mathbb C^{5d+5}$ (They are not global coordinates of $\tilde \mathbb C^{5d+5}$). We may assume
$\tilde c_2$ lies in the neighborhood of the coordinates\footnote{If not, we continue to have successive blow-ups till the preimage of $c_2$ lies 
in the coordinates' neighborhood.}.  So 
$\tilde c_2$ has specific coordinates
\begin{equation}\begin{array}{cc} &
r_0=0 \\
& \theta_i^j=\tilde \theta_i^j, i=0, \cdots, 4,  j=1, \cdots, d, \\
& r_j=y_j\neq 0, j=1, \cdots, 4
\end{array}\end{equation}

Because $t_3, \cdots, t_{5d}$ are distinct and the last $2d$ of them are non-zeros, by the definition 4.2, in this neighborhood of $\tilde \mathbb C^{5d+5}$, we have analytic coordinates
\begin{equation} \theta_0^1, \cdots, \theta_{2}^d, r_0, \cdots, r_3, \xi, \epsilon_1, \cdots, \epsilon_{2d}.
\end{equation}

In section 4, these are for $M$. But here we will use them for the space $U_{\tilde \mathbb C^{5d+5}}$ which can be  analytically identified
with a neighborhood $M$ (See (5.14). This identification will be used again in the second step in a more iterated fashion.). 

\smallskip

Then the finite set 
\begin{equation} \mathcal B=\{dr_l, d\xi, d\theta_i^j, d\epsilon_n \}_{l=0, \cdots, 3, i=0, 1, 2, j=1, \cdots, d, n=1, \cdots, 2d}\end{equation}
 is  a basis for $(T_{\tilde c} U_{\tilde \mathbb C^{5d+5}})^\ast$ where
$\tilde c\in U_{\tilde \mathbb C^{5d+5}}$. 
\bigskip

For this step, it is clear that lemma 5.3 follows from 
\bigskip

\begin{lemma}
Let $\mathcal A$ be the coefficient matrix of $5d+5$ vectors in (5.17) under the basis $\mathcal B$. Then
$\mathcal A$ is non degenerate near $\tilde c_2$( $\mathcal A$ is a $(5d+5)\times (5d+5)$ square matrix).

\end{lemma} 

\bigskip

\begin{proof} of lemma 5.4:  
\par
{\it Set-up}: Assume  the blow-up in (5.13). 
The row vectors in $\mathcal A$ are vectors in (5.17). We place them from top to bottom in the following order
\begin{equation} \begin{array}{cc} &
\pi_1(\psi_3), \\
& \vdots\\
& \pi_1^\ast (\psi_{5d}),\\
& \pi_1^\ast (\psi_{5d+1})\\
& \pi_1^\ast (df_2(c(t_{1}))),\\
& \pi_1^\ast (df_2(c(t_{2}))),\\
& \pi_1^\ast (df_1(c(t_{1}))),\\
& \pi_1^\ast (df_1(c(t_{2}))),\\
& \pi_1^\ast (df_0(c(t_{1}))),\\
& \pi_1^\ast (df_0(c(t_{2}))).
\end{array}\end{equation}
The basis vectors are listed,  from left to right, as 
$$
d\theta_0^2, \cdots, \widehat d\theta_1^1, \cdots,  d\theta_2^d, 
d\epsilon_1, \cdots, d\epsilon_{2d}, dr_0, d\theta_0^1, d\theta_1^1, dr_1, dr_2, dr_3, d\xi
$$
( $\widehat\cdot $  denotes omitting).
Then matrix $\mathcal A$ is the matrix  partitioned as 

\begin{equation}\left(\begin{array}{ccc}
\mathcal A_{1 1} &\mathcal A_{12} &\mathcal A_{13} \\
\mathcal A_{2 1} &\mathcal A_{22} &\mathcal A_{23} \\
\mathcal A_{3 1} &\mathcal A_{32} &\mathcal A_{33}\\
\mathcal A_{4 1} &\mathcal A_{42} &\mathcal A_{43}\\
\mathcal A_{5 1} &\mathcal A_{52} &\mathcal A_{53}
\end{array}\right)\end{equation}
where all $\mathcal A_{ij}$ are matrices of different sizes. We should describe them, one-by-one,  as follows\footnote{As we begin here, only $\mathcal A_{11} $
and $\mathcal A_{13}$ are needed to be specific. Other blocks will be  carefully studied only after (5.47)}:
\par

In the following $ O(1)$ denotes a polynomial function on $U_{\tilde \mathbb C^{5d+5}}$ of
forms 
\begin{equation}
(\theta_i^j-\tilde\theta_i^j) O, (\epsilon_k-\beta_k)O, \ for \ i\leq 2, k\leq 2d
\end{equation}
where $O$ are polynomial functions on $U_{\tilde \mathbb C^{5d+5}}$. In the context, $O$ and
  $O(1)$ could be used for different functions. \par

(I) $ \mathcal A_{11}$: 
it is a $(5d-2)\times (5d-2)$ matrix. Entries are coefficients of 
$$d\theta_i^j, d\epsilon_n, i\leq 2, (i, j)\neq (0, 1), (1, 1), n=1, \cdots, 2d$$ for the vectors 
$$  \pi_1^\ast (\psi_i), i=3, \cdots, 5d.
$$
For this block matrix we use part (a), lemma 5.1 to obtain that 
 its diagonal entries in the order of top-to-bottom are
\begin{equation} 
b_3, \cdots, b_{5d}
\end{equation}
where all $b_i=r_0^2 b_i'$ such that  $b_i'$ are polynomial functions on $U_{\tilde \mathbb C^{5d+5}}$ that do not vanish at $\tilde c_2$.
All the rest of entries are in the form $$r_0^2 O(1),$$ where 
$O(1) $ is defined in (5.31).\par

(II) $\mathcal A_{12}$: it is $(5d-2)\times 1$ matrix. Entries are coefficients of
$dr_0$ for the vectors $$\pi_1^\ast (\psi_i), i=3, \cdots, 5d.$$
So from top-to-bottom they are
$$a_{3}, \cdots, a_{5d}$$
where  $a_l=r_0 O(1)$.

(III) $ \mathcal A_{13}$: it is a $(5d-2)\times 6$ matrix. Entries are coefficients of $$d\theta_0^1, d\theta_1^1, dr_1, dr_2, dr_3, d\xi$$ for the vectors

$$ \pi_1^\ast (\psi_i), i=3, \cdots, 5d.$$
So we use part (b) of lemma 5.1 to obtain that they are all in the form
$$r_0^2 O(1).$$

(IV) $\mathcal A_{21}$: it is a $1\times (5d-2)$ matrix. Entries are coefficients of 
$$d\theta_i^j, d\epsilon_n, i\leq 2, (i, j)\neq (0, 1), (1, 1), n=1, \cdots, 2d$$   for
the vector 
$$\pi_1^\ast (\psi_{5d+1}).$$
So all entries are in the form $$r_0^2 O$$
 where 
$O $ is a polynomial function on $U_{\tilde \mathbb C^{5d+5}}$. 

\par
(V) $\mathcal A_{22}$: it is a $1\times 1$ matrix. It is the coefficient of $dr_0$ for the vector, 
$$\pi_1^\ast (\psi_{5d+1}).$$
It is in the form
$$r_0 O$$ where $O $ is a polynomial function on $U_{\tilde \mathbb C^{5d+5}}$ .
\par

(VI) $\mathcal A_{23}$: it is a $1\times 6$ matrix. Entries are coefficients of 
$$d\theta_0^1, d\theta_1^1, dr_1, dr_2, dr_3, d\xi$$ for the vector,
$$\pi_1^\ast (\psi_{5d+1}).$$

So all entries are in the form $$r_0^2 O$$
 where 
$O $ is a polynomial function on $U_{\tilde \mathbb C^{5d+5}}$. 
\par
(VII) $\mathcal A_{31}$: it is a $2\times (5d-2)$ matrix. Entries are coefficients of
$$d\theta_i^j, d\epsilon_n, i\leq 2, (i, j)\neq (0, 1), (1, 1), n=1, \cdots, 2d$$   for the vectors
$$\pi_1^\ast (df_2(c(t_{1}))),  \pi_1^\ast (df_2(c(t_{2}))).$$
So the entries  are all in the form
$$r_0 O, $$
where 
$O $ is a polynomial function on $U_{\tilde \mathbb C^{5d+5}}$. 

\par

(VIII) $\mathcal A_{32}$: it is $2\times 1$ matrix. Entries are the coefficients of 
$dr_0$ for the vectors
$$\pi_1^\ast (df_2(c(t_{1}))),  \pi_1^\ast (df_2(c(t_{2}))).$$
So they are in
the form of $$ O,$$
where 
$O $ is a polynomial function on $U_{\tilde \mathbb C^{5d+5}}$. 

(IVV) $\mathcal A_{33}$: it is $2\times 6$ matrix. Entries are the coefficients of 
$$d\theta_1^1, d\theta_2^1, dr_1, dr_2, dr_3, d\xi$$ for the vectors
$$\pi_1^\ast (df_2(c(t_{1}))),  \pi_1^\ast (df_2(c(t_{2}))).$$
So all entries are in the forms of 
$$r_0 O,$$
where 
$O $ is a polynomial function on $U_{\tilde \mathbb C^{5d+5}}$. 

(VV) $\mathcal A_{41}$: it is a $2\times (5d-2)$ matrix. Entries are coefficients of
$$d\theta_i^j, d\epsilon_n, i\leq 2, (i, j)\neq (0, 1), (1, 1), n=1, \cdots, 2d$$    for
$$
 \pi_1^\ast (df_1(c(t_{1}))),
 \pi_1^\ast (df_1(c(t_{2}))).
$$
All entries are in the form of $$r_0 O,$$
where 
$O $ is a polynomial function on $U_{\tilde \mathbb C^{5d+5}}$. 
\par

(VVI) $\mathcal A_{42}$: it is a $2\times 1$ matrix. Entries are the coefficients of
$dr_0$ for the vector
$$
 \pi_1^\ast (df_1(c(t_{1}))),
 \pi_1^\ast (df_1(c(t_{2}))).
$$
\par
(VVII) $\mathcal A_{43}$: it is a $2\times 6$ matrix. 
Entries are the coefficients of
$$d\theta_0^1, d\theta_1^1, dr_1, dr_2, dr_3, d\xi$$ for the vector
$$
 \pi_1^\ast (df_1(c(t_{1}))),
 \pi_1^\ast (df_1(c(t_{2}))).
$$
All entries are in the form 
$$r_0 O,$$
where 
$O $ is a polynomial function on $U_{\tilde \mathbb C^{5d+5}}$.

\par

(VVIII) $\mathcal A_{51}$: it is a $2\times (5d-2)$ matrix. Entries are coefficients of
$$d\theta_i^j, d\epsilon_n, i\leq 2, (i, j)\neq (0, 1), (1, 1), n=1, \cdots, 2d$$    for
$$
 \pi_1^\ast (df_0(c(t_{1}))),
 \pi_1^\ast (df_0(c(t_{2}))).
$$
\par

(VVIV) $\mathcal A_{52}$: it is a $2\times 1$ matrix. Entries are the coefficients of
$dr_0$ for the vector
$$
 \pi_1^\ast (df_0(c(t_{1}))),
 \pi_1^\ast (df_0(c(t_{2}))).
$$
\par
(VVV) $\mathcal A_{53}$: it is a $2\times 6$ matrix. 
Entries are the coefficients of
$$d\theta_0^1, d\theta_1^1, dr_1, dr_2, dr_3, d\xi$$ for the vector
$$
 \pi_1^\ast (df_0(c(t_{1}))),
 \pi_1^\ast (df_0(c(t_{2}))).
$$

Next consider the function on $U_{\tilde \mathbb C^{5d+5}}-E_1$ (where $r_0\neq 0$), 
\begin{equation}
\mu(\tilde c)={1\over (r_0)^{10d+1}}|\mathcal A|
\end{equation}
which is the determinant of
the matrix
\begin{equation} \mathcal A_{r_0}=\left(\begin{array}{ccc}
{1\over r_0^2}\mathcal A_{1 1} &{1\over r_0^2}\mathcal A_{12} &{1\over r_0^2}\mathcal A_{13} \\
{1\over r_0}\mathcal A_{2 1} &{1\over r_0}\mathcal A_{22} &{1\over r_0}\mathcal A_{23}\\
{1\over r_0}\mathcal A_{3 1} &{1\over r_0}\mathcal A_{32} &{1\over r_0}\mathcal A_{33}\\
{1\over r_0}\mathcal A_{4 1} &{1\over r_0}\mathcal A_{42} &{1\over r_0}\mathcal A_{43}\\
\mathcal A_{5 1} &\mathcal A_{52} &\mathcal A_{53}\\
\end{array}\right). \end{equation}
To prove lemma 5.4, it is sufficient to prove the determinant 
$$|\mathcal A_{r_0}|\ at\ r_0=0 $$
is non-zero. 
Unfortunately,  the matrix $\mathcal A_{r_0}$ is not a well-defined at $\tilde c_2$ (i.e. $r_0=0$) 
 because some entries involve  ${(t_l-\theta_i^j)\over r_0}$. But those terms do not 
show up in the computation of its determinant $\mu(\tilde c)$. Hence $\mu(\tilde c)$ can be continuously extended to $\tilde c_2$.
  Let's see this. \smallskip

{\it Computation}: 

Notice all entries in 

\begin{equation} \left(\begin{array}{c}
{1\over r_0^2}\mathcal A_{1 1} \\
{1\over r_0}\mathcal A_{2 1} \\
{1\over r_0}\mathcal A_{3 1} \\
{1\over r_0}\mathcal A_{4 1}\\
\mathcal A_{5 1} \\
\end{array}\right). \end{equation}
and 
\begin{equation} {1\over r_0^2}\mathcal A_{13} \end{equation}
can be extended to the entire neighborhood $U_{\tilde \mathbb C^{5d+5}}$, and when evaluated at $\tilde c_2$, 
\begin{equation} {1\over r_0^2}\mathcal A_{1 1}
\end{equation}
is a non-zero diagonal matrix (by part (a), lemma 5.1).  Also notice that 
\begin{equation} {1\over r_0^2}\mathcal A_{13}|_{\tilde c_2}=0. \end{equation}
Therefore in case $r_0\neq 0$, we can perform the row operations on 
the matrix
$$\mathcal A_{r_0}$$ to reduce
\begin{equation} \left(\begin{array}{c}
{1\over r_0}\mathcal A_{2 1} \\
{1\over r_0}\mathcal A_{3 1} \\
{1\over r_0}\mathcal A_{4 1}\\
\mathcal A_{5 1} \\
\end{array}\right) \end{equation}
to zero matrix.
Hence 
\begin{equation}
\mu(\tilde c)={1\over (r_0)^{10d+1}}|\mathcal A|
\end{equation}
is a non-zero multiple of 
\begin{equation} \left|\begin{array}{cc}
{1\over r_0}(\mathcal A_{22}+O_{22}(1)) &{1\over r_0}\mathcal A_{23}\\
{1\over r_0}(\mathcal A_{32}+O_{32}(1)) &{1\over r_0}\mathcal A_{33}\\
{1\over r_0}(\mathcal A_{42}+O_{42}(1)) &{1\over r_0}\mathcal A_{43}\\
\mathcal A_{52}+{1\over r_0}O_{52}(1) &\mathcal A_{53}\\
\end{array}\right|+O(1). \end{equation}
where $O(1), O_{ij}(1)$ represent determinants  of matrices and matrices,  whose entries are local functions vanishing at $\tilde c_2$. Notice that
in (5.41) the block matrices of first row and first column have sizes $1\times 7$ and $7\times 1$ respectively.  
Therefore by the linearity of determinants, 
\begin{equation} \mu(\tilde c)={1\over \rho} \left|\begin{array}{cc}
{1\over r_0}\mathcal A_{22} &{1\over r_0^2}\mathcal A_{23}\\
\mathcal A_{32} &{1\over r_0}\mathcal A_{33}\\
\mathcal A_{42} &{1\over r_0}\mathcal A_{43}\\
r_0\mathcal  A_{52} &\mathcal A_{53}\\
\end{array}\right|+O(1) \end{equation}
for some non-zero complex number $\rho$.    
Now all entries in (5.42) are well-defined functions on $U_{\tilde \mathbb C^{5d+5}}$. 
Therefore it suffices to show the non-degeneracy of this $7\times 7$ matrix 
\begin{equation} \left.\left(\begin{array}{cc}
{1\over r_0}\mathcal A_{22} &{1\over r_0^2}\mathcal A_{23}\\
\mathcal A_{32} &{1\over r_0}\mathcal A_{33}\\
\mathcal A_{42} &{1\over r_0}\mathcal A_{43}\\
r_0\mathcal  A_{52} &\mathcal A_{53}\\
\end{array}\right)\right|_{\tilde c_2}
\end{equation}
at the point $\tilde c_2$. 
Notice the determinant of the matrix (5.43) is also a Jacobian matrix of 7 functions

$$\begin{array}{c} \pi_1^\ast(\psi_{5d+1}), \pi_1^\ast(f_2(c(t_1))), \pi_1^\ast(f_2(c(t_2))),\\
\pi_1^\ast(f_1(c(t_1))),  \pi_1^\ast(f_1(c(t_2))), , \\ 
\pi_1^\ast(f_0(c(t_1))), \pi_1^\ast(f_0(c(t_2)))\end{array} $$
of $7$ variables
$$ \theta_0^1, \theta_1^1, r_0, \cdots, r_3, \xi .$$
By the definition 4.2, 
there is an isomorphic coordinates' system of other $7$ variables
$$ \theta_0^1, \theta_1^1, r_0, \cdots, r_3, r_4.$$
Therefore Jacobians under two different coordinates' system differ only by a non-zero multiple. 
Next we use the same notations $\mathcal A_{ij}$  to denote the coefficient matrices  under the new basis
$$ d\theta_0^1, d\theta_1^1, d r_0, \cdots, d r_3, d r_4.$$
So we switch the matrix  to show the determinant of (5.42) at $\tilde c_2$ under the new basis
\begin{equation}  d\theta_0^1, d\theta_1^1, d r_0, \cdots, d r_3, d r_4 \end{equation}
is non-zero ( replace $\xi$ by $r_4$). \bigskip

Because $t_{5d+1}$ is generic,  also $q$ is generic and $\tilde\theta_i^j$ are distinct, 
the first row vector in (5.42) is generic in $\mathbb C^7$ with respect to other row vectors ( but it is not true under the old basis 
$$ d\theta_0^1, d\theta_1^1, d r_0, \cdots, d r_3, d \xi ).$$
Thus it suffices to prove $6\times 6$ matrix
\begin{equation} \left.\left(\begin{array}{c}
{1\over r_0}\mathcal A_{33}\\
{1\over r_0}\mathcal A_{43}\\
\mathcal A_{53}
\end{array}\right)\right|_{\tilde c_2}
\end{equation}
 is non-degenerate.  

Let $\lambda_1$ be the determinant of (5.45). Also let
$$ g_i(t)=({1\over r_0}\pi_1^\ast ({ z_i\partial  f_1\over \partial z_i}))|_{\tilde c_2(t)}.$$
Using the coordinates in definition 4.2, we compute that this determinant
\begin{equation} \lambda_1=\left. det \left (\begin{array}{c}
{1\over r_0}\mathcal A_{33}\\
{1\over r_0}\mathcal A_{43}\\
\mathcal A_{53}
\end{array}\right) \right|_{\tilde c_2}
\end{equation}
is equal to 

\begin{equation} \lambda_2 \left|\begin{array}{cccccc}
{1\over t_1-\tilde\theta_0^1} &{1\over t_1-\tilde\theta_1^1} &1 &1  & 1 & 1  \\
{1\over t_2-\tilde\theta_0^1} &{1\over t_2-\tilde\theta_1^1} &1 &1 & 1 & 1 \\
{\partial \pi_1^\ast( f_1)(\tilde c_2 (t_1))\over r_0\partial \theta_0^1} &{\partial \pi_1^\ast( f_1)(\tilde c_2 (t_1))\over r_0\partial \theta_1^1} &
g_1(t_1) & g_2(t_1) &
 g_{3}(t_1)& g_4(t_1) \\
{\partial \pi_1^\ast( f_1)(\tilde c_2 (t_2))\over r_0\partial \theta_0^1} &{\partial \pi_1^\ast( f_1)(\tilde c_2 (t_2))\over r_0\partial \theta_1^1} &
g_1(t_2) & g_2(t_2) &
 g_{3}(t_2)& g_4(t_2) \\
 {\partial \pi_1^\ast( f_0)(\tilde c_2 (t_1))\over \partial \theta_0^1} & {\partial \pi_1^\ast( f_0)(\tilde c_2 (t_1))\over \partial \theta_1^1}
& (z_1{\partial f_0 \over \partial z_1})|_{\tilde c_2(t_1)} &(z_2{\partial f_0 \over \partial z_2})|_{\tilde c_2(t_1)} &
 (z_{3}{\partial f_0 \over \partial z_{3}})|_{\tilde c_2(t_1)}& (z_4{\partial f_0 \over \partial z_4})|_{\tilde c_2(t_1)}\\
{\partial \pi_1^\ast( f_0)(\tilde c_2 (t_2))\over \partial \theta_0^1} & {\partial \pi_1^\ast( f_0)(\tilde c_2 (t_2))\over \partial \theta_1^1}
& (z_1{\partial f_0 \over \partial z_1})|_{\tilde c_2(t_2)} &(z_2{\partial f_0 \over \partial z_2})|_{\tilde c_2(t_2)} & 
(z_{3}{\partial f_0 \over \partial z_{3}})|_{\tilde c_2(t_2)}& (z_4{\partial f_0 \over \partial z_4})|_{\tilde c_2(t_2)}
\end{array}\right|. \end{equation}

where 
\begin{equation}
\lambda_2=\prod_{i=1}^2 ({\pi_1^\ast(f_2)\over y_1y_2y_3y_4r_0})|_{\tilde c_2(t_i)}\neq 0. 
\end{equation}

Because $q$ is generic,  the two middle row vectors of (5.47)  

\begin{equation}\begin{array}{cccccc}
({\partial \pi_1^\ast( f_1)(\tilde c_2(t_1))\over r_0\partial \theta_0^1} &{\partial \pi_1^\ast( f_1)(\tilde c_2(t_1))\over r_0\partial \theta_1^1} &
g_1(t_1) & g_2(t_1) &
 g_{3}(t_1)& g_4(t_1)) \\
({\partial \pi_1^\ast( f_1)( \tilde c_2(t_2))\over r_0\partial \theta_0^1} &{\partial \pi_1^\ast( f_1)(\tilde c_2(t_2))\over r_0\partial \theta_1^1} &
g_1(t_2) & g_2(t_2) &
 g_{3}(t_2)& g_4(t_2))
\end{array}\end{equation}
evaluated at $\tilde c_2$ are generic vectors in $\mathbb C^6$ with respect to others. Thus 
to show (5.47) is non-zero, it suffices to show that the minor

\begin{equation} \begin{array}{c} Jac(f_0, c_2)\\
\|\\
\left|\begin{array}{cccccc}
{1\over t_1-\tilde\theta_0^1}  &1  & 1 & 1  \\
{1\over t_2-\tilde\theta_0^1} &1 & 1 & 1 \\
 {\partial \pi_1^\ast( f_0)(\tilde c_2 (t_1))\over \partial \theta_0^1}  & (z_2{\partial f_0 \over \partial z_2})|_{\tilde c_2(t_1)} &
 (z_{3}{\partial f_0 \over \partial z_{3}})|_{\tilde c_2(t_1)}& (z_4{\partial f_0 \over \partial z_4})|_{\tilde c_2(t_1)}\\
{\partial \pi_1^\ast( f_0)(\tilde c_2 (t_2))\over \partial \theta_0^1} & (z_2{\partial f_0 \over \partial z_2})|_{\tilde c_2(t_2)} &
 (z_{3}{\partial f_0 \over \partial z_{3}})|_{\tilde c_2(t_2)}& (z_4{\partial f_0 \over \partial z_4})|_{\tilde c_2(t_2)}
\end{array}\right|\\
\nparallel\\
0. \end{array}\end{equation}
(Remove the 2nd and 3rd columns in (5.47)).
We further compute to have
\begin{equation}\begin{array}{c}  Jac(f_0, c_2)\\
\|\\
({1\over t_1-\tilde\theta_0^1}-{1\over t_2-\tilde\theta_0^1} )\left|\begin{array}{ccc}
1&1&1\\  (z_2{\partial  f_0\over \partial z_2})|_{c_2(t_1)} &
(z_3{\partial  f_0\over \partial z_3})|_{c_2(t_1)} &(z_4{\partial  f_0\over \partial z_4})|_{c_2(t_1)} \\
(z_2{\partial  f_0\over \partial z_2})|_{c_2(t_2)} &
(z_3{\partial  f_0\over \partial z_3})|_{c_2(t_2)} &(z_4{\partial  f_0\over \partial z_4})|_{c_2(t_2)} 
\end{array}\right|. \end{array}
\end{equation}

The determinant $Jac(f_0, c_2)$ is a regular function of $f_0, c_2$ (assuming $t_1, t_2$ are fixed). 
We would like to prove the assertion 

$$  
\left|\begin{array}{ccc}
1&1&1\\  (z_2{\partial  f_0\over \partial z_2})|_{c_2(t_1)} &
(z_3{\partial  f_0\over \partial z_3})|_{c_2(t_1)} &(z_4{\partial  f_0\over \partial z_4})|_{c_2(t_1)} \\
(z_2{\partial  f_0\over \partial z_2})|_{c_2(t_2)} &
(z_3{\partial  f_0\over \partial z_3})|_{c_2(t_2)} &(z_4{\partial  f_0\over \partial z_4})|_{c_2(t_2)} 
\end{array}\right|\neq 0 
$$

for any generic $f_0$ and $c_0$ that has no multiple zeros with coordinates planes $\{z_i=0\}$. 
Let $\Sigma$ be an open subvariety
$$\{c\in M: zeros \ of \ c_i(t)=0\ are \ distinct, i=1, \cdots, 4\}.$$
Consider the family of rational maps

$$V_{f}=\{ c\in \Sigma: Jac(f, c)=0 \}.$$
Notice by the definition $V_{f}$ is a subvariety of $\Sigma$.
Next we consider the fibre $V_{f_1}$ where $f_1$ is the Fermat quintic
$$f_1=z_0^5+\cdots+z_4^5.$$
It is obvious $V_{f_1}$ is empty. Hence $V_{f}$ is empty for generic $f$.  This shows that 
\begin{equation} \left|\begin{array}{ccc}
1&1&1\\  (z_2{\partial  f_0)\over \partial z_2})|_{c_2(t_1)} &
(z_3{\partial  f_0)\over \partial z_3})|_{c_2(t_1)} &(z_4{\partial  f_0)\over \partial z_4})|_{c_2(t_1)} \\
(z_2{\partial  f_0)\over \partial z_2})|_{c_2(t_2)} &
(z_3{\partial  f_0)\over \partial z_3})|_{c_2(t_2)} &(z_4{\partial  f_0)\over \partial z_4})|_{c_2(t_2)} 
\end{array}\right|\neq 0. 
\end{equation}
Therefore 
$$Jac(f_0, c_2)\neq 0.$$
This completes the proof of lemma 5.4, thus the first step of lemma 5.3 in the case 
$c^i_2\neq 0, i\neq 0$. \par

In other cases with more $c_2^i=0$, the proof is identical if we continue the blow-ups along the vanishing components to reach
the $\tilde c_2$ with non-vanishing coordinates $\tilde c_2^i\neq 0, i=0, \cdots, 4$.  The eventual termination of the blow-ups is similar to
that in step 2. So we refer the readers to step 2.

\end{proof}
\bigskip

{\bf Step 2}:  The following formulation of  step 2  is  similar to that of step1. \par
Let's assume $\tilde \theta_i^j$ are not distinct.  
 Notice $\omega(M, \mathbf t)$ depends on the choice of $t_i, i=1, \cdots, 5d+1$.
The main  idea we would like to get across in this paper is that for generic choice of $t_i$, $\omega(M, \mathbf t)$ is non-zero near $c_2$ (excluding $c_2$) .
 However it seems to be completely helpless in the computation of $\omega(M, \mathbf t)$ for generic $t_i$. This situation changes if we use a special choice of $t_i$ such as in step 1, in which case the Jacobian matrix in $\omega(M, \mathbf t)$ breaks down to manageable block matrices in (5.30).  The choice of such $t_i$ requires that
$\tilde \theta_i^j$ are distinct. If they aren't (i.e. $\tilde c_2$ has  multiple intersection points with coordinates' plane), 
 the step 1 will fail because the main block in (5.30) becomes degenerate at all $c\in  P (\Gamma_{\mathbb L})$. 
So to resolve this,  we blow-up these multiple roots, i.e. will use successive blow-ups  to reduce the multiple $\tilde \theta_i^j$ to distinct $\tilde \theta_i^j$. 
They eventually become distinct because this is the case for a generic $c\in P (\Gamma_{\mathbb L})$. 
The pull-back of $\omega(M, \mathbf t)$ by the blow-ups for generic $t_i$ is exactly the same as (5.10) which is that in step 1.

\bigskip

Let's discuss the blow-ups.  They are all analytic which means the total spaces and the centers are analytic in an analytic neighborhood. 
Suppose we already have the blow-up $\tilde \mathbb C^{5d+5}$ as in (5.13). 
 Recall  $\tilde c_2^i\in H^0(\mathcal O_{\mathbf P^1}(d))$ are  sections such that 
$$ div(\tilde c_2^i)=\sum_{j=1}^{d}\tilde \theta_i^j. $$
for $i=0, \cdots, 4$, where $\tilde c_2^0$ is regarded as an element in the nonprojectivzed $H^0(\mathcal O_{\mathbf P^1}(d))$. 
We may assume
multiple zeros are
$$\tilde\theta_\alpha^{\beta}$$
among all $\tilde\theta_i^j$, 
where $\alpha, \beta$ are finite numbers less than $5$ and $d+1$ respectively. Let 
$$m(\tilde c_2)$$ be the largest number of the pairs $(\alpha, \beta)$ such that $\tilde \theta_\alpha^\beta$ are the same.  The worst case is 
$m(\tilde c_2)=5d$, in which case $\tilde c_2$ represents a constant map $\mathbf P^1\to \mathbf P^4$. This multiplicity is well-defined for the map $\tilde c_2$, but depends on the coordinates $z_i$. 
Next we define the successive blow-ups which are meant to reduce the multiplicity $m(\tilde c_2)$. \par
{\it First blow-up}.  First we should note, as in definition 4.1, the coordinates $r_i, \theta_i^j$ are well-defined analytic coordinates for $\tilde M$ around the point $\tilde c_2$, even if $\theta_i^j$ are not distinct.\footnote{ This is actually true in a Zariski topology. The analytic blow-ups we are using here can be replaced by 
the algebraic ones.} 
Let $\Delta_{d+1}$ be an open disk of $\mathbb C^{d+1}$ (not around the origin). We let 
$$\Delta^{(0)}\simeq \Delta_{d+1}\times\cdots\times \Delta_{d+1}$$
be an open set of $\tilde M$ around $\tilde c_2$. Each copy $\Delta_{d+1}$ is an open set
of $H^0(\mathcal O_{\mathbf P^1}(d))$.  Thus each point of $\Delta^{(0)}$ is a five tuple of 
sections in $H^0(\mathcal O_{\mathbf P^1}(d))$. We use the same coordinates $r_i, \theta_i^j$ 
for each copy of $\Delta^{(0)}$. So a point in  $\Delta^{(0)}$ can be expressed as
$$\begin{array}{c} (\tilde c^0(t), \cdots, \tilde c^4(t))\\
where\ \tilde c^i(t)=r_i\sum_{j=0}^d (t-\theta_i^j)^j, r_i\neq 0.\end{array} $$
We define  $B_2$ to be the subvariety of $\Delta^{(0)}$, 
 consisting of all points of $\Delta^{(0)}$ that has at least $m(\tilde c_2)$ common zeros $\theta_i^j\in \mathbf P^1$.
Note $B_2$ is a union of planes of dimension $$5(d+1)-m(\tilde c_2)+1,$$
and 
$$\tilde c_2$$ is a point on $B_2$ which could be a singular point lying on the intersection of multiple planes. 
But for the simplicity we may assume it is a smooth point that lies on a single plane.   

We blow-up $\Delta^{(0)}$ along $B_2$ to obtain the first blow-up map
$$\begin{array}{ccc}
N & \stackrel{\pi_2} \rightarrow &  \Delta^{(0)}.
\end{array}$$

Let $c_2^{(1)}$ be the inverse  of $\tilde c_2$ in the strict transform of $\tilde P(\Gamma_{\mathbb L_1})$.
Using the coordinates $\theta_i^j, r_i$ (in definition 4.1) above,  we can find an open set $$\Delta^{(1)}\subset N$$ centered around $ c_2^{(1)}$ such that
 $\Delta^{(1)}$  is analytically isomorphic to
\begin{equation}
\Delta^{(1)} \simeq  \Delta_{d+1}\times \cdots \times \Delta_{d+1},\end{equation}
 where each copy $\mathbb C^{d+1}$ corresponds to $\theta_i^j$ with the same $i$. 
The following is the detailed description of $\Delta^{(1)}$.
In general the isomorphism (5.53) is not unique. But in (5.53), we would like to choose a particular isomorphism 
such that  outside of $\pi_2^{-1}(B_2)$, 
$$\pi_2^\ast(r_i), \pi_2^\ast(\theta_i^j)$$ are the analytic coordinates of $i$-th copy of 
$$\Delta_{d+1}\times \cdots \times \Delta_{d+1}.$$
We should denote $$\pi_2^\ast(r_i), \pi_2^\ast(\theta_i^j)$$ 
by $$r_i, (\theta_i^j)^{(1)}.$$
( $(\theta_i^j)^{(1)}$ is just the strict transform of $\pi_2^\ast(\theta_i^j)$ ). 
This implies that 
$$\pi_2=\zeta_0^1\times \zeta_1^1\times\cdots\times \zeta_4^1$$
outside of the center, where each $\zeta_i^1$ is an analytic automorphism of $\Delta_{d+1}$ (though there is a relation among all $\zeta_i$ ). 
Next we fix the embedding  on each copy $\Delta_{d+1}$, 
\begin{equation}\begin{array}{ccc} \Delta_{d+1} &  \hookrightarrow  & H^0(\mathcal O_{\mathbf P^1}(d))\\
(r_i, (\theta_i^1)^{(1)}, \cdots, (\theta_i^d)^{(1)})& \rightarrow &  r_i \sum_{j=1}^d(t-(\theta_i^j)^{(1)})= (c^i)^{(1)}(t). \end{array}\end{equation}

Then we obtain the new evaluation map $Z^{(1)}$, 
\begin{equation}\begin{array}{ccc}
\Delta^{(1)}\times \mathbf P^1 &\stackrel{Z^{(1)}} \rightarrow & \mathbf P^4\\
((r_i, \theta_i^j), t) &\rightarrow & [( c^0)^{(1)}(t), \cdots, ( c^4)^{(1)}(t)].
\end{array}\end{equation}

The map $Z^{(1)}$  turns $\Delta^{(1)}$ into a family of maps
\footnote{ Because the non isomorphism on the center, rational curves in this family do not lie on hypersurfaces in the family $\mathbb L$. So they are  not all from the family 
$\tilde P(\Gamma_{\mathbb L})$.  But $\tilde c_2^{(1)}$ is a specialization of the generic member in this family $\Delta^{(1)}$ of rational curves. 
More geometrically, $\Delta^{(1)}$ consists of a collection of some open paths in $\tilde P(\Gamma_{\mathbb L})$ whose open location is filled up with all partial derivatives of the rational curves along these paths.}, $$\mathbf P^1\to \mathbf P^4.$$

{\it Continuing blow-ups}.  If  $ c_2^{(1)}(t)$ still has multiple zeros with the planes $\{z_i=0\}$, 
i.e. $(  c_2^i)^{(1)}(t)=0$ for all $i$ have multiple zeros or equivalently 
$$m( c_2^{(1)})>1, $$ we continue the blow-ups along the multiple zeros as above. 
This type of  blow-ups can continue. 
The last blow-up is obtained when the multiplicity $m( c_2^{(\kappa)})$ is reduced to $1$ ( need to repeat (5.53) and (5.54) in each blow-up to have
the well-defined $m(c_2^{(\kappa)})$ ). 
To see that 
 $m(\tilde c_2)$ will be reduced to 1 (i.e. the blow-ups will stop), we go back to the first blow-up. 
Let $\tilde P(\Gamma_{\mathbb L_1})^{(1)}$ be the strict transform of $\tilde P(\Gamma_{\mathbb L_1})$ under the first blow-up, and
$ c_2^{(1)}$ be a chosen inverse of $\tilde c_2$ in $\tilde P(\Gamma_{\mathbb L_1})^{(1)}$.  Because there is a $c$ in 
$\tilde P(\Gamma_{\mathbb L_1})$, but not in $B_2$ (this is the birationality of $c_0$),  $\tilde P(\Gamma_{\mathbb L_1})$
does not lie in the  $B_2$. Thus the exceptional divisor $\mathcal D$ of $\tilde P(\Gamma_{\mathbb L_1})^{(1)}$ does not lie in 
	the tangent bundle of $B_2$ ( $\mathcal D$ is a subvariety in the projectivization of the normal bundle of $B_2$),
where the projectivization of the
normal bundle
of $B_2$ exactly corresponds to the collection of those new rational curves $c^{(1)}$ at infinity (If $\tilde c_2$ is not a smooth point of $B_2$, we just do finite such blow-ups to strictly reduce $m(\tilde c_2)$).
Therefore  $c_2^{(1)}$ when regarded as a map  $\mathbf P^1\to \mathbf P^4$  has multiplicity
$$m( c_2^{(1)})$$ strictly  less than $m(\tilde c_2)$.
After such successive blow-ups, we obtain the inverse $c_2^{(\kappa)}$ of $\tilde c_2$ whose multiplicity
$$m( c_2^{(\kappa)})$$ is one. \footnote {If generic $c\in P(\Gamma_{\mathbb L})$ is a multiple-cover map, such successive blow-ups may not be terminated because  $\tilde P(\Gamma_{\mathbb L_1})\subset B_2$ in the first step. This is one of main reasons why $c_0$ can not be a multiple-cover map.}
We let $\pi_3$ be the composition of all such blow-ups. 
So we obtain the birational map

$$ \begin{array}{ccc}
\Delta^{(\kappa)} &\stackrel{\pi_3}\rightarrow & \Delta^0.
\end{array}$$
where  $\Delta^{(\kappa)}$ is an open set of  the final blow-up, centered at $c_2^{(\kappa)}$, and
a determined  analytic isomorphism 
\begin{equation}
\Delta^{(\kappa)}\simeq \Delta_{d+1}\times \cdots\times\Delta_{d+1} .
\end{equation}
We denote the natural affine coordinates of $$\Delta^{(\kappa)}$$
in (5.56)  by $c^{(\kappa)}=( (c_i^j)^{(\kappa)})_{i=0, j=1}^{i=4, j=d}$. The formula (5.55) can be written as
$$ \begin{array}{ccc}
\Delta^{(\kappa)} \times \mathbf P^1 &\rightarrow & \mathbf P^4\\
( (c_i^j)^{(\kappa)})_{i=0, j=1}^{i=4, j=d} &\rightarrow & \biggl( (c^0)^{(\kappa)}(t), \cdots,  (c^4)^{(\kappa)}(t)\biggr).
\end{array} $$

\bigskip
{\it Computation of the pull-back by blow-ups}. 
Next we apply the successive blow-ups $\pi_3$ to prove the lemma. 
We would like to show that after the blow-ups we have the same set-up as for the step 1. 
Let $\tilde P(\Gamma_{\mathbb L})^{(\kappa)}$ be the strict transform of $\tilde P(\Gamma_{\mathbb L})$ under $\pi_3$. 
  Then $\Delta^{(\kappa)}\cap \tilde P(\Gamma_{\mathbb L})^{(\kappa)}$ parametrizes a family of rational curves in $\mathbf P^{4}$ 
\footnote {These rational curves outside of exceptional locus lie on the same hypersurfaces from $\mathbb L$, but will be specialized to ``higher partial derivatives" of the rational curves.  }.

These successive blow-ups can be stated in the following diagram

\begin{equation}\begin{array}{ccccccc}
\tilde P(\Gamma_{\mathbb L})^{(\kappa)}& \stackrel{P} \leftarrow & W &\subset  &\Delta^{(\kappa)}\times \mathbb L &\stackrel{Proj.}  \rightarrow & \mathbb L\\
\downarrow\scriptstyle{ \pi_3} & & \downarrow&& \downarrow\scriptstyle{ \pi_3\times {identity}}  && ||\\
\tilde P(\Gamma_{\mathbb L})  &\stackrel{P} \leftarrow & \tilde \Gamma_{\mathbb L} &\subset  & \Delta^{(0)}\times \mathbb L  & \stackrel{Proj.}  \rightarrow & \mathbb L
\end{array}\end{equation}
where $W$ is the strict transform of $\tilde \Gamma_{\mathbb L}$ and $\pi_3$ is birational.  
Notice that outside of $\pi_3^{-1}(B_2)$, $\pi_3$ can be decomposed as
an isomorphism
$$\zeta_0^{\kappa}\times \zeta_1^{\kappa}\times\cdots\times \zeta_4^{\kappa}$$
where each $\zeta_i^{\kappa}$ is an analytic automorphism of $\Delta^{(\kappa)}$.

Next we repeat the construction in step 1 to construct the differential form for  $\tilde P(\Gamma_{\mathbb L})^{(\kappa)}$. 
Let $q$ be a generic homogeneous quadratic polynomial in $z_0, \cdots, z_4$ as before. 
Let
$$
f_3=z_0z_1 z_2( \delta_1 q+\delta_2 z_{3} z_4),
$$
and  $t_1, t_2$ be generic numbers satisfying
$$
\left|\begin{array}{cc} q|_{ c_2^{(\kappa)}(t_1)}  &  ( c_2^3)^{(\kappa)}(t_1) ( c_2^4)^{(\kappa)}(t_1)\\
q|_{ c_2^{(\kappa)} (t_2)}& ( c_2^3)^{(\kappa)}(t_2) ( c_2^4)^{(\kappa)} (t_2))\end{array}\right|=0.
$$
where $c_2^{(\kappa)}=(  c_2^0)^{(\kappa)}, \cdots,  ( c_2^4)^{(\kappa)})$, and
$$\delta_1=\left|\begin{array}{cc} f_0( c_2^{(\kappa)}(t_1))  &   f_2( c_2^{(\kappa)}(t_1))\\
f_0 ( c_2^{(\kappa)} (t_2))& f_2( c_2^{(\kappa)}(t_2))\end{array}\right|, \quad \delta_2=\left|\begin{array}{cc} f_1 ( c_2^{(\kappa)}(t_1))  &   f_0( c_2^{(\kappa)}(t_1))\\
f_1 ( c_2^{(\kappa)}(t_2))& f_0(  c_2^{(\kappa)})(t_2))\end{array}\right|.$$
Let $ (\tilde\theta_i^j)^{(\kappa)}$ be the zeros of 
$( c_2^i)^{(\kappa)}(t)=0, i=0, \cdots, 4$. 
Then let 
$t_3, \cdots, t_{5d}$ be zeros of
\begin{equation}\begin{array}{c} f_3(z)|_{ c_2^{(\kappa)}(t)}\\
\|\\
( c_2^0)^{(\kappa)} (t) \cdots( c_2^{2})^{(\kappa)}(t) \biggl ( 
\delta_1
q|_{ c_2^{(\kappa)}(t)}+  \delta_2 ( c_2^{3})^{(\kappa)}(t) ( c_2^4)^{(\kappa)}(t)\biggr)\\
\|\\
0.\end{array}
\end{equation}
other than $ (\tilde\theta_0^1)^{(\kappa)},  (\tilde\theta_1^1)^{(\kappa)}$.  Let $t_{5d+1}$ be generic.  
Notice that the only difference between the set-up in step 1 and that in this step is that
rational curves  $( c)^{(\kappa)}$, i.e. in this step we have replaced 
each $ c=( c_0, \cdots,  c_4)\in \tilde M$ in step 1 by $ c^{(\kappa)} $ in step 2. 
( The only difference between $c$ and $ c^{(\kappa)}$ lies in their specializations which are outside their centers. This means that the interested rational curves $ c_2^{(\kappa)}$ in the variables $ c^{(\kappa)}$ do not exist in variables $c$). 
With such a special choice of $t_i$, 
we can repeat the step 1 for  the same quintics $f_i(z)$ to construct
$ \omega^{(\kappa)}(\Delta^{(\kappa)}, \mathbf t)$ where $\omega^{(\kappa)}(\Delta^{(\kappa)}, \mathbf t)$ has the current choice of $t_i$ and the variables are $c^{(\kappa) }$ for $\Delta^{(\kappa)}$. Specifically, we define
\begin{equation} \phi_i^{(\kappa)}= d
\left|  \begin{array}{ccc} f_2( c^{(\kappa)}(t_i)) & f_1(c^{(\kappa)}(t_i)) & f_0( c^{(\kappa)}(t_i))\\
f_2( c^{(\kappa)}(t_1)) & f_1( c^{(\kappa)}(t_1)) & f_0( c^{(\kappa)}(t_1))\\
f_2( c^{(\kappa)}(t_2)) & f_1( c^{(\kappa)}(t_2)) & f_0( c^{(\kappa)}(t_2))
\end{array}\right|\end{equation}
for $i=3, \cdots, 5d+1$
\begin{equation}
\omega^{(\kappa)}(\Delta^{(\kappa)}, \mathbf t)=\wedge_{i=3}^{5d+1} \phi_i^{(\kappa)}\in H^0(\Omega_{\Delta^{(\kappa)}}).
\end{equation}
Then we repeat the same process (but without blow-up $\pi_1$) in step 1 to have the Jacobian 
matrix 
$$\mathcal A^{(\kappa)}$$ 
with respect to $\omega^{(\kappa)}(\Delta^{(\kappa)}, \mathbf t)$ and the point $ c_2^{(\kappa)}$.
Then since $ c^{(\kappa)}$ are specialized to rational curves that have non multiple zeros with coordinates' planes, we use the same process in step 1 to calculate $det( \mathcal A^{(\kappa)})$. We obtain 
same result as in lemma 5.4,  i.e.
\begin{equation}
det( \mathcal A^{(\kappa)})\neq 0. 
\end{equation}
at the center point $ c_2^{(\kappa)}$. 
Next we discuss the relation between 
$det(\mathcal A)$ and $det( \mathcal A^{(\kappa)}) $.
We would like see that,  as in step 1
$det(\mathcal A)$ is for the Jacobian at the generic points of $\tilde P(\Gamma_L)$, 
 $det( \mathcal A^{(\kappa)}) $ now is for the Jacobian at the generic points of $( \tilde P(\Gamma_L))^{(\kappa)}$.
Notice that $\pi_3$ is an isomorphism outside of exceptional divisor $\pi_3^{-1}(B_2)$. \par
We claim that  the equations
$$\left|  \begin{array}{ccc} f_2( c^{(\kappa)}(t_i)) & f_1( c^{(\kappa)}(t_i)) & f_0( c^{(\kappa)}(t_i))\\
f_2( c^{(\kappa)}(t_1)) & f_1( c^{(\kappa)}(t_1)) & f_0( c^{(\kappa)}(t_1))\\
f_2( c^{(\kappa)}(t_2)) & f_1( c^{(\kappa)}(t_2)) & f_0( c^{(\kappa)}(t_2))
\end{array}\right|=0$$
$i=3, \cdots, 5d+1$ defines the scheme 
$$\tilde P(\Gamma_{\mathbb L})^{(\kappa)}=\pi_3^{-1} (\tilde P(\Gamma_{\mathbb L}))$$
in $\Delta^{(\kappa)}-\pi_3^{-1}(B_2)$. 
To prove the claim, we notice for a polynomial $F$ in $\mathbb C^5$ (with coordinates $\tilde c$), there is a well-defined function
$$\pi_3^\ast (F(\tilde c(t))$$
on $\Delta^{(\kappa)}$ (with coordinates $ c^{(\kappa)}$), 
which is exactly 
$$F( c^{(\kappa)})$$ outside of the center of $\pi_3$.  Outside of center of $\pi_3$, this can be  expressed as 
 
$$\pi_3^{-1} ( \{\tilde c: F(\tilde c(t))=0\})=\{ c^{(\kappa)}: F( c^{(\kappa)})=0\}.$$
Apply the quintics to $F$.  
This isomorphism $\pi_3$ shows that 
  on the open set 
\begin{equation} \Delta^{(\kappa)}\cap \tilde P(\Gamma_{\mathbb L})^{(\kappa)}-\pi_3^{-1}(B_2) \end{equation}

\begin{equation} (\pi_3)^\ast (det(\mathcal A))=g \cdot det( \mathcal A^{(\kappa)}) \end{equation}
where $g$ is a function nowhere zero on 
$$ \tilde P(\Gamma_{\mathbb L})^{(\kappa)}-\pi_3^{-1}(B_2).$$
Because of (5.61), $det( \mathcal A)\neq 0$ at a point in 
$$ \pi_1\circ \pi_3(\Delta^{(\kappa)})\cap P(\Gamma_{\mathbb L})-\{c_2\}. $$

This proves lemma 5.3.

\end{proof}

\bigskip

\begin{proof}  of lemma 5.2: 
We first calculate 1-form $\phi_i$ on $M$ evaluated at general $c_g\in M$. 
For $i=3, \cdots, 5d+1$, 
\begin{equation}\begin{array}{cc}
  \phi_i= d
\left|  \begin{array}{ccc} f_2(c(t_i)) & f_1(c(t_i)) & f_0(c(t_i))\\
f_2(c(t_1)) & f_1(c(t_1)) & f_0(c(t_1))\\
f_2(c(t_2)) & f_1(c(t_2)) & f_0(c(t_2))
 \end{array}\right|&\\ 
=\left|  \begin{array}{cc} f_0(c_g(t_1)) &  f_2(c_g(t_1))\\
f_0(c_g(t_2)) & f_2(c_g(t_2))\end{array}\right| df_1(c(t_i))+\left|  \begin{array}{cc} f_2(c_g(t_1)) &  f_1(c_g(t_1))\\
f_2(c_g(t_2)) & f_1(c_g(t_2))\end{array}\right| df_0(c(t_i))& \\  +\left|  \begin{array}{cc} f_1(c_g(t_1)) &  f_0(c_g(t_1))\\
f_1(c_g(t_2)) & f_0(c_g(t_2))\end{array}\right| df_2(c(t_i))
+ \sum_{l=0, j=1}^{l=2, j=2} h_{lj}^i(c_g) df_l(c(t_j)) & \end{array}\end{equation}

where $h_{lj}^i$ are polynomials in $c$. \par

The lemma 5.3 implies that 
for a generic choice of $$f_0, f_1, f_2, c_g, t_1, \cdots, t_{5d+1},$$

\begin{equation}\begin{array}{cc} & 
\left|  \begin{array}{cc} f_0(c_g(t_1)) &  f_2(c_g(t_1))\\
f_0(c_g(t_2)) & f_2(c_g(t_2))\end{array}\right| df_1(c(t_i))+\left|  \begin{array}{cc} f_2(c_g(t_1)) &  f_1(c_g(t_1))\\
f_2(c_g(t_2)) & f_1(c_g(t_2))\end{array}\right| df_0(c(t_i))\\ &
+\left|  \begin{array}{cc} f_1(c_g(t_1)) &  f_0(c_g(t_1))\\
f_1(c_g(t_2)) & f_0(c_g(t_2))\end{array}\right| df_2(c(t_i)) \end{array}\end{equation} for
$i=3, \cdots, 5d+1$, and
\begin{equation} df_l(c(t_j)), l=0, 1, 2, j=1, 2
\end{equation}
are $5d+5$ linearly independent vectors in $(T_{c_g}M)^\ast$ for a generic $c_g\in P(\Gamma_{\mathbb L})$ (not at special point $c_2$ ).
 i.e. they
form a basis of the vector space $(T_{c_g}M)^\ast$.

This implies the set of 1-forms $\{\phi_i\}_{ i=3, \cdots, 5d+1}$ is a linearly independent set in $(T_{c_g}M)^\ast$ for
generic $c_g\in P(\Gamma_{\mathbb L})$.  Thus $\omega(M, \mathbf t)$ is nowhere zero when
it is evaluated at  generic points of $P(\Gamma_{\mathbb L})$.  The lemma 5.2 is proved.

\end{proof}

\bigskip

\subsection{Ranks of differential sheaves}
\quad\smallskip

\bigskip

\begin{proof} of proposition 1.4: 
Let $\mathcal N$ be the submodule of global sections, $H^0(\Omega_{M})$
generated by
elements
\begin{equation} \phi_i= d
\left|  \begin{array}{ccc} f_2(c(t_1)) & f_1(c(t_1)) & f_0(c(t_1))\\
f_2(c(t_2)) & f_1(c(t_2)) & f_0(c(t_2))\\
f_2(c(t_i)) & f_1(c(t_i)) & f_0(c(t_i))\end{array}\right|\end{equation}
for $i=3, \cdots, 5d+1$. 
Recall that 
$$\left|  \begin{array}{ccc} f_2(c(t_1)) & f_1(c(t_1)) & f_0(c(t_1))\\
f_2(c(t_2)) & f_1(c(t_2)) & f_0(c(t_2))\\
f_2(c(t_i)) & f_1(c(t_i)) & f_0(c(t_i))\end{array}\right|=0,$$
for $i=3, \cdots, 5d+1$
define the scheme $P(\Gamma_{\mathbb L})$ for a small $\mathbb L$.  By   proposition 8.12 in [7], II, \begin{equation}
\widetilde {({H^0(\Omega_{M})\over \mathcal N})}\otimes \mathcal O_{P(\Gamma_{\mathbb L})}
\simeq \Omega_{P(\Gamma_{\mathbb L})},
\end{equation}
where $\widetilde {(\cdot)}$ denotes the sheaf associated to the module $(\cdot)$.

Therefore

\begin{equation}\begin{array}{cc} &
 {({H^0(\Omega_{M})\otimes k(c_g)\over \mathcal N\otimes k(c_g)})}
\simeq \Omega_{P(\Gamma_{\mathbb L})}\otimes k(c_g)\\&
=(\Omega_{P(\Gamma_{\mathbb L})})|_{(\{c_g\})},\end{array}\end{equation}
where $k(c_g)=\mathbb C$ is the residue field at generic $$c_g\in P(\Gamma_{\mathbb L}).$$ 
Notice two sides of (5.69) are finitely dimensional  linear spaces over $\mathbb C$.
\begin{equation}\begin{array}{cc} &
dim_{\mathbb C} ((\Omega_{P(\Gamma_{\mathbb L})})|_{(\{c_g\})})\\
&=dim_{\mathbb C}(H^0(\Omega_{M})\otimes k(c_g))
- dim(\mathcal N\otimes k(c_g))\end{array}
\end{equation}

Since \begin{equation}
dim_{\mathbb C} ((\Omega_{P(\Gamma_{\mathbb L})})|_{(\{c_g\})}))=dim(T_{c_g}P(\Gamma_{\mathbb L}))
\end{equation}

\begin{equation}
dim(T_{c_g}P(\Gamma_{\mathbb L}))=dim(M) -dim(\mathcal N\otimes k(c_g)).
\end{equation} 

By lemma 5.2, 
$$dim(\mathcal N\otimes k(c_g))=deg(\omega(M, \mathbf t))=5d-1.$$

The proposition 1.4 is proved.

\end{proof}

\bigskip

 \begin{proof} of theorem 1.1. By lemma 5.2, 
$$dim(\mathcal N\otimes k(c_g))=5d-1.$$

Thus by proposition 1.4, 
$$dim(T_{c_g}P(\Gamma_{\mathbb L}))=5d+5-(5d-1)=6.$$

Then using lemma 3.2 part (c) and lemma 3.3, we obtain 
$$H^1(N_{c_0/X_0})=0,$$
for generic $c_0\in Hom(\mathbf P^1, X_0)$. Next we show 
that theorem 1.1 does not need $c_0$ to be generic because there are only finitely many rational curves of degree $d$ on each 
generic $X_0$. More specifically
by lemma 3.3 again,  if $dim (P(\Gamma_{\mathbb L_1}))=5$, 
$\Gamma_{f_0}$ has to have dimension $4$,  which must be the union of orbits 
of $GL(2)$. Thus each component of  $Hom_{bir}(\mathbf P^1, X_0)$ is
an orbit of $GL(2)$. This implies $$H^1(N_{c_0/X_0})=0,$$ for all
$c_0\in Hom(\mathbf P^1, X_0)$. The theorem is proved.
\end{proof}

\section{Examples\\
--Vainsencher's and Chen's rational curves}

{\bf Example 6.1} (Vainsencher's rational curves)
 \par
This example provides an evidence to theorem 1.1.
In [8], Vainsencher constructed irreducible, degree 5, nodal curves $C_0$ on a generic quintic $f_0$ by taking plane sections of the quintic. 
Let $c_0$ be its normalization. By
our theorem 1.1, $c_0$ is an immersion and 
\begin{equation}
N_{c_0/X_0}\simeq \mathcal O_{\mathbf P^1}(-1)\oplus \mathcal O_{\mathbf P^1}(-1).
\end{equation}

Indeed these were proved by Cox and Katz in [4], by using a different method. Their method is based on Clemens' deformation idea.   Their understanding of $c_0$ on $f_0$  was achieved by a concrete construction 
of special $c_0$, $f_0$ and by using a computer program for the
last verification of the $26\times 30$ matrix.  It is easy to check that the rational maps $c_0$ they constructed are immersions.  \par
Furthermore
our result shows 
\begin{equation} dim(T_{c_0}\Gamma_{f_0})=4.
\end{equation}
Because of the equation (6.1), $C_0$ can't deform in $f_0$. Thus
$\Gamma_{f_0}$ consists of multiple orbits isomorphic to $GL(2)(c_0)$.   Theorem 1.1 also shows that there will not be any scheme-theoretical multiplicity associated to the orbits. However the number of these orbits is not accessible because the degree of each orbit in $\mathbf P(M)$
could be different. This number is related to  Gromov-Witten invariants.

\bigskip

{\bf Example 6.2} (Chen's rational curves) \par

This is an example on $K$-3 surfaces. 
In [1], Chen constructed  nodal rational curves $C_0$ of degree $4d$ for each natural number $d$,  that lie on the 
generic hypersurfaces $f_0$ of degree $4$ in $\mathbf P^3$ ($f_0$ is a K-3 surface).   At first we may have an impression that this is against our intuition. 
Because it is similar to rational curves on generic  quintic threefolds that we can have naive counting: on a generic quartic  hypersurface $f_0$ of $\mathbf P^3$ , there will be $4d+1$ conditions imposed the rational curves on $f_0$, while the dimension of the moduli space of rational curves in $\mathbf P^3$ (modulo $PGL(2)$ action) is only $4d$. Thus the naive counting concludes that there will not be any rational curves on $f_0$.  But it was proved by Mori, Mukai, etc.,  and  Chen ([1]) that rational curves on $f_0$ exist and they are all nodal. Our proof is closely related to this counting, and our construction of $\omega(M, \mathbf t)$ can be carried out in $\mathbf P^3$ for Chen's case.  
But theorem 1.1 does not hold because proposition 1.3 fails.  This failure is not expected by the naive dimension count, but it is a reminder of a fact that  the  generic quartics are not generic in the moduli space of complex  structures.    \par
 
Chen's  construction has a similar
flavor of Vainsencher's rational curves  above. They were obtained by taking hyperplane sections of $K$-3 surfaces. 
Intrinsically Vainsecher's and Chen's rational curves look similar.  For instance they are all plane sections, and are all immersed,  nodal rational curves.   So what invariant distinguishes one from the other?  Section 5 shows that this invariant may not be the invariant of the intrinsic rational curves, it  addresses the structure of the moduli space of rational curves for underlined families of varieties. 
 More specifically, it is deduced from the differential form $\omega(M, \mathbf t)$ (defined in (1.8) ). The $\omega$ itself is not a moduli invariant, but the 
zero locus $\{\omega(M, \mathbf t)=0\}$ is, and furthermore $\{\omega(M, \mathbf t)=0\}$ is independent of generic $t_i, i=1, \cdots, 5d+1$. In Chen's situation, $\omega(M, \mathbf t)$ turns out to be  identically zero on  $P(\Gamma_{\mathbb L})$,  but in Vainsencher's  it is not.   Beyond Chen's cases, it is not clear that which homology classes of rational curves would have or would not have vanishing $\omega(M, \mathbf t)$.

\end{document}